\documentclass[letterpaper, paper,12pt]{AAS}		
\usepackage[colorlinks=true, pdfstartview=FitV, linkcolor=black, citecolor= black, urlcolor= black]{hyperref}
\usepackage{bm}
\usepackage{amsmath}
\usepackage{amsthm}
\usepackage{subfigure}
\usepackage{graphicx}
\usepackage{titlesec}
\usepackage{mathtools}

\usepackage{overcite}
\usepackage{footnpag}			      	
\usepackage{amsfonts}
\usepackage{amssymb}
\usepackage{ragged2e}
\usepackage{rotating}
\usepackage{cleveref}
\usepackage{multicol}
\newcommand{\V}[1]{\ensuremath{\textbf{\text{#1}}}} 
\newcommand{\tr}[1]{{#1}^{\ensuremath{\mathsf{T}}}} 
\usepackage{xcolor}

\PaperNumber{24-480}

\usepackage{float}
\usepackage{comment}

\newtheorem{theorem}{Theorem}
\newtheorem{lemma}{Lemma}

\newtheorem{assumption}{Assumption}
\newtheorem{corollary}{Corollary}

\begin{document}

\title{Optimal Control with Lyapunov Stability Guarantees for Space Applications}

\author{Abhijeet\thanks{Department of Aerospace Engineering,
Texas A \& M University,
710 Ross St,
College Station, 77843, Texas, USA},
Mohamed Naveed Gul Mohamed\footnotemark[1],
Aayushman Sharma\footnotemark[1],
Suman Chakravorty\footnotemark[1]
}

\maketitle{}

\begin{abstract}

This paper investigates the infinite horizon optimal control problem (OCP) for space applications characterized by nonlinear dynamics. The proposed approach divides the problem into a finite horizon OCP with a regularized terminal cost, guiding the system towards a terminal set, and an infinite horizon linear regulation phase within this set. This strategy guarantees global asymptotic stability under specific assumptions. Our method maintains the system's fully nonlinear dynamics until it reaches the terminal set, where the system dynamics is linearized. As the terminal set converges to the origin, the difference in optimal cost incurred reduces to zero, guaranteeing an efficient and stable solution. The approach is tested through simulations on three problems: spacecraft attitude control, rendezvous maneuver, and soft landing. In spacecraft attitude control, we focus on achieving precise orientation and stabilization. For rendezvous maneuvers, we address the navigation of a chaser to meet a target spacecraft. For the soft landing problem, we ensure a controlled descent and touchdown on a planetary surface. We provide numerical results confirming the effectiveness of the proposed method in managing these nonlinear dynamics problems, offering robust solutions essential for successful space missions.
\end{abstract}

\section{Introduction}
\justifying
The successful exploration and utilization of space requires advanced control strategies to ensure the success of various missions. From the precise orientation of spacecraft to the delicate maneuvers required for docking and landing, control technique plays a pivotal role in overcoming the inherent challenges of the space environment. While traditional control methods can be effective in certain scenarios, they often fall short in providing an optimal and stable feedback solution needed for long-duration missions within the dynamic and unpredictable space environment \cite{blackmore2013lossless, SMC_att, ipopt_attitude}. In addition, achieving global asymptotic stability is particularly valuable for space missions, where the system has to reach a specific terminal state irrespective of its initial state. \cite{mohamed2024optimalsolutioninfinitehorizon}.

Over time, various methods and strategies have been devised to tackle control problems in space applications. The shooting method, for instance, has become notable for solving optimal control problems and is particularly effective for applications like attitude control \cite{att_con_Shoot}, trajectory optimization \cite{traj_opt_shoot}, and others \cite{shoot_gen, shoot_reentry}. Despite its benefits, the shooting method has limitations, such as poor flexibility in adapting to disturbances or deviations from expected states and high sensitivity to initial guesses \cite{shooting_drawback,rao2009survey}. As an alternative, direct methods like Sequential Quadratic Programming (SQP) or Interior Point methods \cite{Num_opt} are frequently used. These techniques are advantageous because they accommodate constraints, making them useful for space missions. They find applications in soft-landing maneuvers \cite{blackmore2013lossless}, trajectory optimization \cite{SQP_traj_opt}, rendezvous and docking operations \cite{rend_IP} \emph{etc}\cite{ipopt_attitude}. Nevertheless, these techniques fail to ensure global asymptotic stability (GAS) or provide feedback for trajectory adjustments, which is crucial in space-related applications. Any deviation in the trajectory could lead to disastrous consequences when using such methods \cite{drawback_sqp}. To overcome these shortcomings, reinforcement learning (RL) has gained traction as a promising solution, utilizing learning to refine and optimize control strategies through environmental interactions \cite{GAUDET_2020, deep_RL_space}. However, RL approaches are highly data-dependent, and their outcomes can be inconsistent\cite{wang2022searchfeedbackreinforcementlearning}.

Considering the limitations of current methods for solving OCPs and the need for long-term stability, we reformulate the problem as an infinite horizon OCP. To solve this, we adopt a tractable approach \cite{mohamed2024optimalsolutioninfinitehorizon,mohamedCDC2023} that ensures feedback and guarantees global asymptotic stability. Our approach uses the Iterative Linear Quadratic Regulator (iLQR), an indirect method, to tackle the optimal control problem while incorporating feedback for system stabilization \cite{ILQG_tassa2012synthesis}. However, due to the infinite horizon nature of the problem, iLQR cannot be deployed directly. To address this, we divide the solution into two stages by introducing a transfer time that transitions the problem from a nonlinear OCP to an LQR regulation problem. We optimize the control inputs as well as the transfer time to guarantee a smooth transition to the terminal set, where LQR regulation is applied. Initially, iLQR guides the infinite-horizon problem to a terminal set, followed by stabilization with the LQR controller. The free-final-time formulation for nonlinear OCP satisfies the Bellman equation, offering a Control Lyapunov Function (CLF) and thus, ensuring global asymptotic stability. By using iLQR, we obtain feedback in the nonlinear phase, while solving the LQR problem in the terminal set ensures feedback at the final stage. We use the Linear Quadratic Regulator (LQR) for regulation for its reliability in linear systems \cite{bryson}. The nonlinear dynamics of space systems, while challenging, are effectively managed by iLQR, which iteratively adjusts control inputs by solving a sequence of locally linearized quadratic problems \cite{ilqg_todorov}. LQR then provides feedback and ensures a stable solution in the terminal linear phase \cite{mine}. It should be noted that the transfer time to the terminal set is also optimized in the algorithm. This is done to ensure minimal cost difference between the nonlinear OCP and the algorithm deployed in the paper, also ensuring a smooth transfer \cite{mohamedCDC2023}. Moreover, the infinite horizon cost acts as a Lyapunov function, which implies that the system can be driven toward the goal from any state. By ensuring this stability, our algorithm offers a solution that accommodates uncertainties, enhancing the reliability and safety of space missions. The assumptions, proofs, and algorithm description are discussed in detail in various sections of the paper.

To demonstrate the effectiveness of our proposed algorithm, we apply it to three critical space applications: spacecraft attitude control, rendezvous maneuver, and soft landing. These problems encompass a wide range of control challenges, from maintaining precise orientation and stability in attitude control to navigating and docking with another spacecraft in rendezvous maneuvers, and ensuring a controlled descent and touchdown in soft landing operations. Each problem requires a control solution capable of managing nonlinear dynamics and ensuring stable, optimal performance throughout the mission.

The remainder of this article is organized as follows: Section \ref{sec:background} provides an in-depth background on the algorithm, including its theoretical foundations and practical considerations.  We also discuss stability analysis and criteria for selecting the transition point to ensure minimal cost difference and effective linearization. Section \ref{sec:methodology} delves into the mathematical formulation of the iLQR method, its iterative process, and the conditions for transitioning to LQR control. Following this, Section \ref{sec:sys_dyn} outlines the specific dynamics and control requirements for the three case studies: spacecraft attitude control, rendezvous maneuver, and soft landing. We present the mathematical models for each problem, highlighting the nonlinear dynamics and control objectives. Next, we provide numerical results, showcasing simulation outcomes and the performance of our algorithm across different scenarios. The results demonstrate the effectiveness of our approach in managing nonlinear dynamics and achieving global asymptotic stability, validating its potential for practical space applications.

\section{Background}\label{sec:background}
We model the problem as an infinite-horizon optimal control problem. This encapsulates free-final time problems typically used in aerospace applications. The theory used to model the problems discussed in this work below is inspired from previous work on infinite horizon nonlinear control \cite{mohamedCDC2023, mohamed2024optimalsolutioninfinitehorizon}. 
 We summarise it below for completeness.

\subsection{Problem Formulation}
The infinite-horizon optimal control problem can be written in the form
\begin{subequations}\label{eq.IHOCP}
\begin{align}
    J_{\infty}^*(\V{x}) = \min_{\{\V{u}_t\}} \sum_{t=0}^{\infty} c(\V{x}_t, \V{u}_t);~~ &\text{given} ~\V{x}_0 = \V{x} \label{eq.infinite_horizon_cost_func} \\
    \text{subject to the dynamics: } \V{x}_{t+1} &= f(\V{x}_t, \V{u}_t), \label{eq.dynamics}
\end{align}
\end{subequations}
where $\V{x}_t \in \mathbb{R}^n$ represents the state vector of the dynamical system, $\V{u}_t \in \mathbb{R}^p$ represents the control input to the dynamical system, and $c(\V{x}_t, \V{u}_t)$ is the incremental cost. The problem described is an infinite horizon optimal control problem (IH-OCP), which is intractable due to the inherent property of dealing with an infinite time horizon. 

It is well-known that the infinite horizon cost-to-go, $J_{\infty}(\cdot)$, satisfies bellman equation if there exists a solution to the IH-OCP \eqref{eq.IHOCP}  \cite[Ch.7]{bertsekas_vol1}:
\begin{align}
    J_{\infty}^*(\V{x}) = \min_{\V{u}} \{ c(\V{x}, \V{u}) + J_{\infty}^*(f(\V{x}, \V{u})) \}. \label{eq.bellman}
\end{align}

Given that $c(\V{x},\V{u}) > 0 ~\forall ~(\V{x},\V{u})\neq (0,0)$, $J^*_{\infty}(\cdot)$ acts as a Control Lyapunov Function (CLF) for the dynamical system \eqref{eq.dynamics}, and thus, the control feedback policy implicitly defined by the optimal cost-to-go function $J^*_{\infty}(\cdot)$, globally asymptotically stabilizes the dynamical system \eqref{eq.dynamics}.  The proof is simple and demonstrated below.

\begin{figure}[t]
    \centering
    \includegraphics[width = \linewidth]{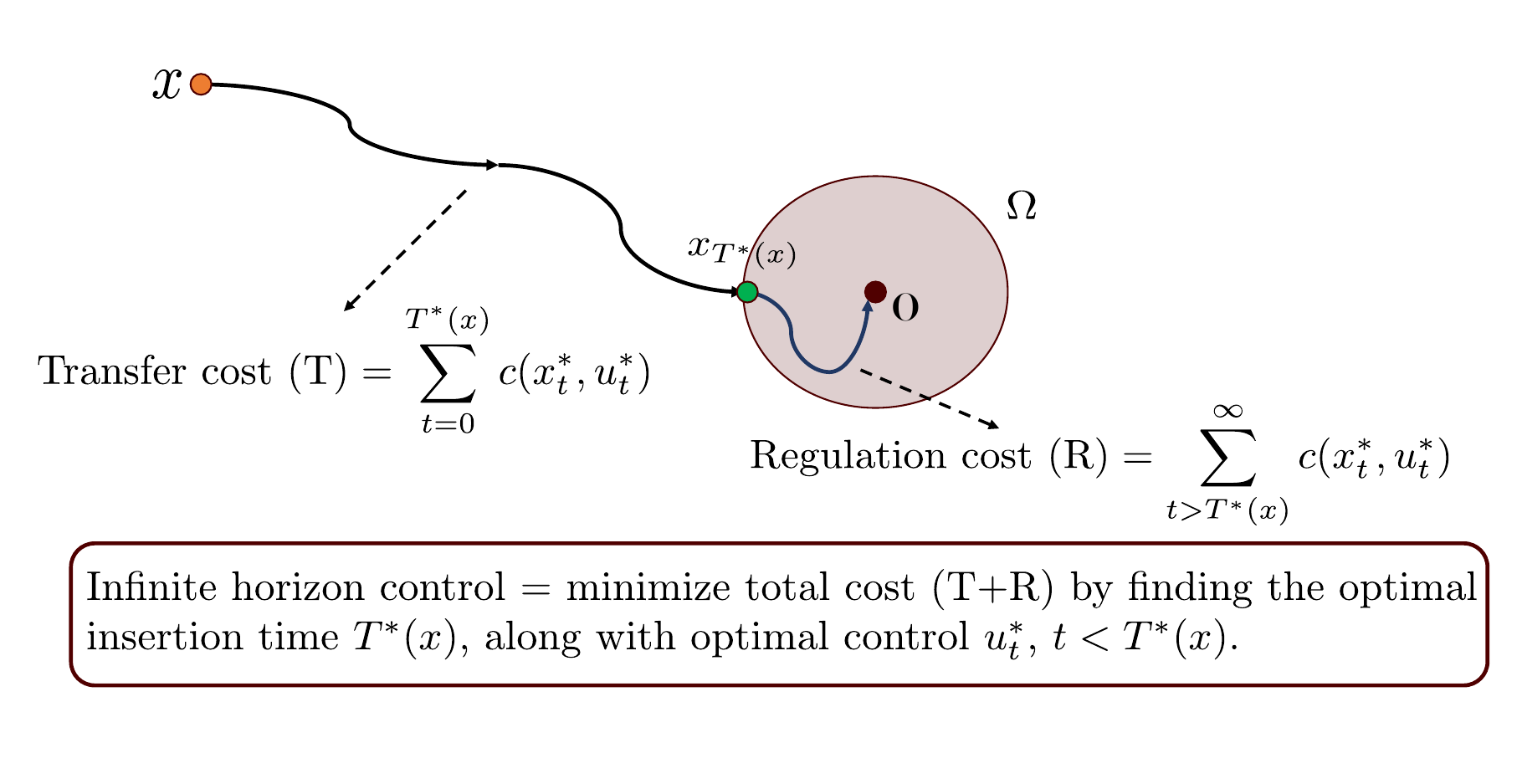}
    \caption{Schematic illustrating the strategy to solve the infinite horizon optimal control problem}
    \label{inf_horizon}
\end{figure}

\begin{corollary} \label{corollary.2}
    Let $J^*_\infty(\V{x})$ satisfy the Bellman equation \eqref{eq.bellman}, then it is a control Lyapunov function for the system in \eqref{eq.dynamics} that renders the origin globally asymptotically stable.
\end{corollary}

Further, suppose that if there exists a $J_{\infty}(\cdot)$ such that it satisfies the Bellman Equation (not necessarily optimal)
\begin{align}
    J_{\infty}(\V{x}) = \min_{\V{u}} \{ c(\V{x}, \V{u}) + J_{\infty}(f(\V{x}, \V{u})) \}, \label{eq.suboptimal_bellman}
\end{align}
then $J_{\infty}(\cdot)$ also is a CLF that renders the origin globally asymptotically stable (GAS).

Thus, our goal for this work is to develop a tractable approach to solving \eqref{eq.IHOCP} by transforming it into a finite horizon problem. Another goal for us in solving \eqref{eq.IHOCP} is to construct CLFs as in \eqref{eq.bellman}/ \eqref{eq.suboptimal_bellman}, such that they render the origin GAS.

We shall make the following assumptions for the rest of the paper. 
\begin{assumption}\label{assump.1 cost}
We assume that the cost function $c(\V{x},\V{u})$ has a global minimum at $(\V{x},\V{u}) = (0,0)$, \textit{i.e.}, $\frac{\partial c}{\partial \V{x}} \Bigr|_{x=0, u=0} = 0$ and $\frac{\partial c}{\partial \V{u}}\Bigr|_{x=0, u=0} = 0$, $c(0,0) = 0$, and $c(\V{x},\V{u}) > 0$ $\forall ~(\V{x},\V{u})\neq(0,0)$.
\end{assumption}
\begin{assumption}\label{assump.2 controllability}
    We assume that given any $\V{x}\in \mathbb{R}^n$, and any $\Omega \subset \mathbb{R}^n$, such that the origin is in $\Omega$, $\exists$ a control sequence $\{\V{u}_t\}_{t=0}^{T(x)}$,  that ensures $\V{x}_{T(x)} \in \Omega$ for some $T(x) < \infty$, under the dynamics defined above (\eqref{eq.dynamics}). 
\end{assumption}
Assumption~\ref{assump.2 controllability} is a controllability assumption that ensures that any state can be controlled into entering the region $\Omega$ in finite time. 
\begin{assumption}\label{assump.3 linear controllability}
    We assume that the linearization of the dynamical system \eqref{eq.dynamics} around $(0,0)$, is controllable.
\end{assumption}

Given assumptions~\ref{assump.1 cost} and \ref{assump.3 linear controllability}, we can define the optimal ``linear" infinite-horizon problem:
\begin{subequations}\label{eq.terminal_lqr}
\begin{align}
    \bar{J}_\infty(\V{x}) &= \min_{\{\V{u}_t \}} \sum_{t=0}^{\infty} (\tr{\V{x}}_t Q \V{x}_t + \tr{\V{u}}_t R \V{u}_t), \label{eq.quad_cost} \\
    \text{subject to:} ~& \V{x}_{t+1} = A \V{x}_t + B \V{u}_t, 
\end{align}
\end{subequations}
where, $(Q,R)$ and $(A,B)$ are obtained by performing a quadratic expansion of $c(\V{x},\V{u})$, and a linear expansion of the dynamics in \eqref{eq.dynamics} around the origin $(\V{x},\V{u}) = (0,0)$.

Note that owing to the linear controllability assumption \ref{assump.3 linear controllability}, $\bar{J}_\infty(\cdot)$ above may be found by solving the stationary algebraic Riccati equation, resulting in 
$
     \bar{J}_\infty(\V{x}) = \tr{\V{x}} P_{\infty} \V{x},
$
where $ P_{\infty} $ is the solution of the stationary Riccati equation.

\subsection{Solution to the Infinite Horizon Optimal Control Problem.}\label{section:sol_IHOCP}
We now define a finite horizon construction to IH-OCP that will use the first hitting time to the set $\Omega_M$, where $\Omega_M = \{\V{x} ~|~ \bar{J}_\infty(\V{x}) \leq M \}$, as the time horizon and whose cost will satisfy the Bellman equation. The construction is suboptimal to the IH-OCP, but we show that the cost of this new construction converges to the true IH-OCP cost in the limit $M \rightarrow 0$.
We call this the alternate construction optimal control problem (AC-OCP), and it is defined as:
\begin{align}
    J^M_\infty(\V{x}) &= \min_{ \{\V{u}_t\}_{t=0}^{T-1} , T} \left( \sum_{t=0}^{T-1} c(\V{x}_t, \V{u}_t) + \max(\bar{J}_\infty(\V{x}_T), M) \right) \label{eq:ACOCP}\tag{AC-OCP} \\
    \text{subject to:} &~\V{x}_{t+1} = f(\V{x}_t, \V{u}_t),\nonumber \\
    & \V{x}_T \in \Omega_M,  ~\text{and given}~ \V{x}_0 = \V{x}.  \nonumber
\end{align}

The above problem has a free final time $T$ that needs to be optimized over in conjunction
with the control actions. The free final time will prove crucial to showing the cost function is a CLF and it converges to the optimal IH-OCP cost. We show the following result.

\begin{theorem}\label{theorem:cost convergence}
    The AC-OCP cost $ J^M_\infty(\V{x})$ converges to the IH-OCP cost $J^*_{\infty}(\V{x})$ in the limit $M \rightarrow 0$, \textit{i.e.}, $$ \lim_{M\rightarrow 0 } J^M_\infty(\V{x}) = J^*_{\infty}(\V{x}),$$ assuming that $J^*_\infty(\cdot)$ is continuous at the origin.
\end{theorem}
The following results also holds.
\begin{lemma}\label{lemma:Bellman equation}
    The cost-to-go of AC-OCP $J^M_\infty(\V{x})$ satisfies the Bellman equation for all initial states $x \notin \Omega_M$, and hence renders the origin globally asymptotic stable.
\end{lemma}


In a nutshell, by converting the infinite horizon control problem into a finite horizon framework, we can manage and optimize the system within a practical and computationally feasible scope. During the first phase, we extend the problem over a long horizon and analyze the behavior of the system under different transfer times. This phase is crucial as it allows us to identify an optimal transfer time that minimizes the total cost without being excessively affected by the regulation cost. The second phase involves ensuring that the set 
$\Omega_{M}$. $\Omega_{M}$ is sufficiently near the system's equilibrium point. This proximity is vital because it enables the linearization of the system dynamics around the equilibrium, simplifying the control design while maintaining accuracy. By ensuring that the regulation cost remains small, we achieve a smooth transition from the nonlinear regime to a linear one. This method leverages the benefits of both finite and infinite horizon approaches. It provides a practical solution framework that is computationally efficient and theoretically robust. By ensuring that the terminal state $x_{T}$ is near equilibrium, we not only achieve global asymptotic stability but also enhance the system's performance and reliability over extended operational periods. This innovative approach bridges the gap between theoretical optimality and practical applicability, making it highly suitable for complex space applications where long-term stability and performance are paramount.  

In some cases, Assumption~\ref{assump.3 linear controllability} is violated. For example, the linearization for nonholonomic systems is not controllable, and some systems do not satisfy the LTI assumption near the final state. For such cases, we assume that the terminal set is forward invariant, as defined below. For these types of systems, similar results as given in Theorem~\ref{theorem:cost convergence} and Lemma~\ref{lemma:Bellman equation} hold\cite{mohamed2024optimalsolutioninfinitehorizon}. 

\begin{assumption}\label{assump.4 forward invariance}
    There exists a control policy $\pi(\cdot): \mathbb{R}^n \rightarrow \mathbb{R}^p$ that makes the set $\Omega_M$ forward invariant under the dynamics in \eqref{eq.dynamics}, \textit{i.e.},
    $f(x,\pi(x)) \in \Omega_M, ~\forall~x\in\Omega_M$. Also, let $c(x,\pi(x)) \eqcolon c^{\pi}(x) \leq \delta ~\forall~ x\in\Omega_M.$ Further $c(x,u) > \delta, ~\forall~ x \notin \Omega_M.$ Here, $\delta$ is a function of $M$, \textit{i.e.}, $\delta = \delta(M)$.
\end{assumption}

\section{Methodology}\label{sec:methodology}
We solve the problem in \ref{eq:ACOCP} indirectly by solving a finite horizon OCP \eqref{eq:FHOCP} and sweeping through different time horizons till the constraint $\V{x}_T \in \Omega_M$ is satisfied. So, we start with a small time $T$ and increase till it satisfies the constraint. The smallest horizon $T$ that satisfies the constraint is the first hitting time for the set $\Omega_M$. To recall, the set $\Omega_M$ is the terminal set where the linear controller is optimal. To check if $\V{x}_T$ is in the set $\Omega_M$, we compare the expected cost-to-go $\bar{J}_{\infty}(\V{x}_T)$ and the actual cost incurred by applying the linear controller on the system and see if they are within a threshold.  
\begin{align}
    J^T_\infty(\V{x}) &= \min_{ \{\V{u}_t\}_{t=0}^{T-1} } \sum_{t=0}^{T-1} c(\V{x}_t, \V{u}_t) + \bar{J}_\infty(\V{x}_T)\label{eq:FHOCP}\\
    \text{subject to:} &~\V{x}_{t+1} = f(\V{x}_t, \V{u}_t),~\text{given}~ \V{x}_0 = \V{x}.  \nonumber
\end{align}
To solve the finite horizon OCP \eqref{eq:FHOCP}, we use the iterative Linear Quadratic Regulator approach (iLQR) \cite{ILQG_tassa2012synthesis}. ILQR is an iterative optimization technique used to solve nonlinear OCPs. It starts with an initial guess for the control sequence $\{\V{u}_0, \V{u}_1, \cdots, \V{u}_{T-1} \}$ and iteratively improves it. It computes the neighboring extremal by using the quadratic expansion of the cost and linear expansion of the dynamics around a trajectory and solves a linear quadratic problem to compute the gains needed to update the control sequence. ILQR is shown to be equivalent to Sequential Quadratic Programming (SQP) \cite{wang2022searchfeedbackreinforcementlearning} in principle but outperforms SQP in computational efficiency owing to the recursive structure in the backward pass computation as opposed to the block computation done in SQP.

Given the optimal control problem in \eqref{eq:FHOCP}, and taking $c(\V{x}_t,\V{u}_t) = l(\V{x}_t) + \frac{1}{2}\tr{\V{u}}_t R \V{u}_t$,  where $l_t(\cdot)$ is the incremental state cost function. The iLQR algorithm consists of the following steps:

\noindent
1. \textbf{Forward Pass.}
Given the nominal trajectory at the previous iteration- $\V{u}^{k}_{0:T-1}$, iLQR gains - $\{k_{0:T-1}, K_{0:T-1}\}$, line search parameter $\alpha$.\\
  Start from $t = 0$, $\bar{\V{x}}_0 = \V{x}_0$.\\
    \begin{align}
        \V{u}^{k+1}_{t} &= \V{u}_t^{k} + \alpha k_t + K_{t} ({\V{x}^{k+1}_t} - {\V{x}_t^{k}}),\\
        {\V{x}^{k+1}_{t+1}} &= f(\V{x}^{k+1}_{t}, {\V{u}_t^{k+1}}),
    \end{align}
\noindent    
2. \textbf{Backward Pass.}
Let $J(\V{x}_T) = \bar{J}_\infty(\V{x}_T)$, 
  compute $J_{\V{x}_T}$ and $J_{\V{x}_T \V{x}_T}$ using terminal conditions. Perform the following steps backward in time for $t=\{T-1, \cdots, 0\}$. 
  First, we compute the partials of the $Q$ function. The $Q(\delta \V{x}_t, \delta \V{u}_t)$ function is the state control value function of the neighboring optimal control problem around a trajectory $(\V{x}_t, \V{u}_t)$.
   \begin{align*}
    Q_{\V{x}_t} &= l_{\V{x}_t}  + \tr{f}_{\V{x}_t} J_{\V{x}_{t+1}},\\
    Q_{\V{u}_t} &= R \V{u}_t + \tr{f}_{\V{u}_t} J_{\V{x}_{t+1}},\\
    Q_{\V{x}_t \V{x}_t} &= l_{\V{x}_t \V{x}_t} + \tr{f}_{\V{x}_t} J_{\V{x}_{t+1} \V{x}_{t+1}} f_{\V{x}_t}, \\
    Q_{\V{u}_t \V{x}_t} &= \tr{f}_{\V{u}_t} (J_{\V{x}_{t+1} \V{x}_{t+1}}) f_{\V{x}_t}, \\
    Q_{\V{u}_t \V{u}_t} &= R + \tr{f}_{\V{u}_t} (J_{\V{x}_{t+1} \V{x}_{t+1}}) f_{\V{u}_t}.
    \end{align*}
Computing the iLQR gains which is used for updating the control,     
    \begin{align}
        k_t &= -Q_{\V{u}_t \V{u}_t}^{-1} Q_{\V{u}_t},\\
        K_t &= -Q_{\V{u}_t \V{u}_t}^{-1} Q_{\V{u}_t \V{x}_t}.
    \end{align}
Compute the partials for the cost-to-go for the previous time-step,    
  \begin{align}
          J_{\V{x}_t} &= Q_{\V{x}_t} + \tr{K}_{t} Q_{\V{u}_t \V{u}_t} k_t + \tr{K}_t Q_{\V{u}_t} + \tr{Q}_{\V{u}_t \V{x}_t} k_t,\\
          J_{\V{x}_t \V{x}_t} &= Q_{\V{x}_t \V{x}_t} + \tr{K}_t Q_{\V{u}_t \V{u}_t} K_t + \tr{K}_t Q_{\V{u}_t \V{x}_t} + \tr{Q}_{\V{u}_t \V{x}_t} K_t.
     \end{align}

\section{System Dynamics}\label{sec:sys_dyn}
The algorithm described above is highly applicable to a variety of space missions. To demonstrate its effectiveness, we have selected three distinct problems as case studies. The following subsections provide a detailed description of the model and system dynamics for each of these problems. It is important to note that while the dynamics are initially formulated as continuous-time ordinary differential equations, for the purposes of this study, we utilize the discrete system dynamics as expressed by equation \eqref{eq.dynamics}.

In each of these case studies, we explore different aspects of space mission challenges, showcasing how the algorithm can be applied to real-world scenarios. The problems have been chosen to illustrate the versatility and robustness of the algorithm in handling complex dynamics and control tasks. Each subsection will delve into specific problem formulations.

\subsection{Attitude Control}
The first problem we consider is the attitude control of a spacecraft in low Earth orbit. For this scenario, we neglect the effect of any disturbance torque. The state vector for this problem consists of six components: [$\psi$, $\theta$, $\phi$, $\omega_{1}$, $\omega_{2}$, $\omega_{3}$]$^{T}$. The governing dynamical equations are:
\begin{subequations}
\begin{align}
\begin{bmatrix}
    \dot{\psi}\\
    \dot{\theta}\\
    \dot{\phi}
\end{bmatrix} =& \frac{1}{\cos{\theta}}\begin{bmatrix}
    0 && \sin{\phi} && \cos{\phi}\\
    0 && \cos{\theta}\cos{\phi} && -\cos{\theta}\sin{\phi}\\
    \cos{\theta} && \sin{\theta}\sin{\phi} && \sin{\theta}\cos{\phi}
\end{bmatrix}\begin{bmatrix}
    \omega_{1}\\
    \omega_{2}\\
    \omega_{3}
\end{bmatrix},\\
\boldsymbol{\dot{\omega}} =& -J^{-1}(\boldsymbol{\omega} \times J \boldsymbol{\omega}) + J^{-1}\boldsymbol{M},
\end{align}
    \label{eq:att_con_dyn}
\end{subequations}

where [$\psi$, $\theta$, $\phi$]$^{T}$ represents 3-2-1 Euler rotation angles. $ \boldsymbol{\omega} $ = $\begin{bmatrix}
     \omega_{1} & \omega_{2} & \omega_{3}
\end{bmatrix} \in \mathbb{R}^3 $ denotes the angular velocity vector in body frame. $J$ is the moment of inertia of the spacecraft, and $\V{M} \in \mathbb{R}^3$ is a vector of control inputs. 

\subsection{Rendezvous Maneuver}
We also consider the rendezvous of a chaser spacecraft with a target spacecraft. The spacecraft model is fully actuated, meaning it has thrust along all three axes. Both the target and the chaser are in elliptical orbits. The state equations are as follows:
\begin{subequations}
    \begin{align}
        \boldsymbol{\dot{e}_{r}} =& \boldsymbol{v_{t}} - \boldsymbol{v_{c}},\\
        \boldsymbol{\dot{e}_{v}} =& -\frac{\mu}{R_{t}^{3}}\boldsymbol{r_{t}} + \frac{\mu}{R_{c}^{3}}\boldsymbol{r_{c}} - \frac{\boldsymbol{u}}{m},\\
        \dot{m} =& -\alpha ||\boldsymbol{u}||,
    \end{align} 
    \label{eq:rend}
\end{subequations}
where $\boldsymbol{e_{r}} = \boldsymbol{r_{t}} - \boldsymbol{r_{c}}$ is the relative error in distance of the chaser with respect to target spacecraft. $\boldsymbol{e_{v}} = \boldsymbol{v_{t}} - \boldsymbol{v_{c}}$ is the relative error in velocity of the chaser with respect to the target spacecraft. $\boldsymbol{r_{t}}$, $\boldsymbol{v_{t}}$ represent the position and velocity vector of the target in the earth-centered inertial frame. $\boldsymbol{r_{c}}$, $\boldsymbol{v_{c}}$ are the position and velocity vector of the chaser in earth-centered inertial frame. $R_{t}$ and $R_{c}$ are the distance of chaser and target from the center of the earth. $\boldsymbol{u} $ =  [$u_{1}$, $u_{2}$, $u_{3}$]$^{T}$ $\in \mathbb{R}^{3}$ denote the control input, and $\mu$ = $398600$ is a constant.

However, the seven equations provided are insufficient to fully propagate the dynamics of the chaser and target spacecraft. Additional equations are needed to provide enough information for this propagation, as the existing equations are expressed in terms of the position and velocity of both the chaser and target in the Earth-centered inertial frame. Therefore, we propagate the position and velocity of the target spacecraft using the following dynamics:
\begin{subequations}
    \begin{align}
        \boldsymbol{\dot{r_{t}}} =& \boldsymbol{v}_{t},\\
        \boldsymbol{\dot{v_{t}}} =& -\frac{\mu}{R_{t}^{3}}\boldsymbol{r_{t}}.
    \end{align}
\end{subequations}
\subsection{Soft-Landing}
The final problem we address is the soft landing of a lander on Mars. The governing dynamical equations are:

\begin{subequations}
    \begin{align}
        \begin{bmatrix}
            \dot{\psi}\\
            \dot{\theta}\\
            \dot{\phi}
        \end{bmatrix} =& \frac{1}{\cos{\theta}}\begin{bmatrix}
            0 && \sin{\phi} && \cos{\phi}\\
            0 && \cos{\theta}\cos{\phi} && -\cos{\theta}\sin{\phi}\\
            \cos{\theta} && \sin{\theta}\sin{\phi} && \sin{\theta}\cos{\phi}
        \end{bmatrix}\begin{bmatrix}
            \omega_{1}\\
            \omega_{2}\\
            \omega_{3}
        \end{bmatrix},\\
        \boldsymbol{\dot{\omega}} =& -J^{-1}(\boldsymbol{\omega} \times J \boldsymbol{\omega}) + J^{-1}\boldsymbol{M},\\
        \boldsymbol{\dot{r}} =& \boldsymbol{v},\\
        \boldsymbol{\dot{v}} =& \frac{\boldsymbol{u}}{m}
         + \boldsymbol{g}_{ref},\\
        \dot{m} =& -\frac{||\boldsymbol{u}||}{I_{sp}g_{ref}},
    \end{align}
\end{subequations}
where [$\psi$, $\theta$, $\phi$]$^{T}$ represents 3-2-1 Euler rotation angles, $\boldsymbol{\omega}$ = [$\omega_{1}$, $\omega_{2}$, $\omega_{3}$]$\in \mathbb{R}^{3}$ is the angular velocity vector [$\omega_{1}$, $\omega_{2}$, $\omega_{3}$], $J$ is the moment of inertia of the spacecraft, $\boldsymbol{M} \in \mathbb{R}^{3}$ is a vector of moment control inputs, and $\boldsymbol{u}$ is a vector of control inputs. $\boldsymbol{r}$ and $\boldsymbol{v}$ are the position and velocity vectors of the lander with respect to an inertial frame with origin as the landing point. $I_{sp}$ is the specific impulse. $\boldsymbol{g}_{ref}$ is the gravity on mars given by $[0,0,-3.7114]^{T}$.

\section{Numerical Results}
The proposed theory is simulated for the dynamical systems described in the above section. The simulations were done in MATLAB, using Euler integration to propagate the non-linear dynamics of the described systems. The initial and terminal conditions are mentioned in Table \ref{tab:init_cond}. The simulation conditions, and the results are discussed in detail in the subsequent subsections. The incremental cost is assumed to be quadratic for all applications, \textit{i.e.}, $c(x,u)$ = $\frac{1}{2}(x^{T}Qx + u^{T}Ru)$.

\begin{table}[!htbp]
    \centering
    \begin{tabular}{|c|c|c|}
    \hline
        System & Initial State & Goal State \\
        \hline
        Attitude control & [$85.94^{o}$, $-68.75^{o}$, $-120.32^{o}$, $5.72^{o}/s$, $-5.72^{o}/s$, $2.86^{o}/s$] & $\boldsymbol{O}_{6 \times 1}$ \\
        \hline
        Soft-Landing &[$22.91^{o}$, $17.18^{o}$, $11.45^{o}$, $5.72^{o}/s$, $11.45^{o}/s$, $-11.45^{o}/s$] & $\boldsymbol{O}_{6 \times 1}$\\
        (attitude params.) && \\
        \hline
        Soft-Landing &[$300 m$, $-200 m$, $1000 m$, $100 m/s$, $120 m/s$, $0 m/s$] & $\boldsymbol{O}_{6 \times 1}$\\
        (position \& velocity) &&\\
        \hline
    \end{tabular}
    \caption{Initial state and goal state for attitude control and soft-landing problem.}
    \label{tab:init_cond}
\end{table}
\subsection{Attitude Control}
The system dynamics have been comprehensively detailed in the preceding section (Eq. \eqref{eq:att_con_dyn}). In this problem, the objective is to maneuver a spacecraft from an initial state to a final state as mentioned in Table \ref{tab:init_cond}. The moment of inertia matrix is given as $J=\text{diag}[4500,2000,7500]$. The maneuver is designed to be executed over a horizon of 200 seconds, with the system dynamics discretized using a step size of 0.1 seconds.

\begin{figure}[!htbp]
    \centering
    \subfigure[Total cost incurred as a function of transfer time.]{\includegraphics[width=0.48\textwidth]{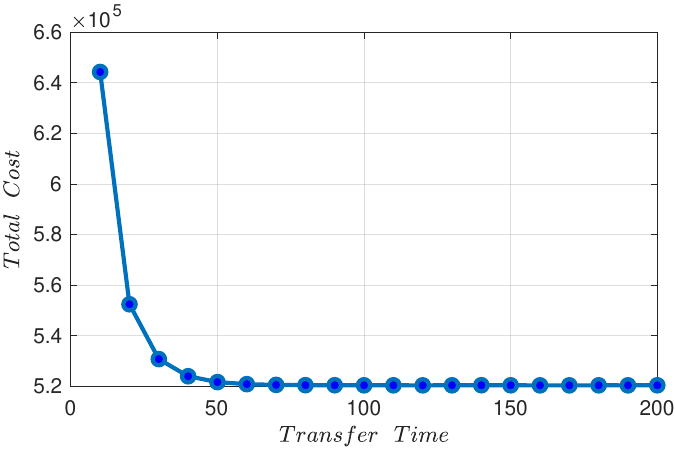}} 
    \subfigure[Cost incurred by linearized system as a function of transfer time.]{\includegraphics[width=0.48\textwidth]{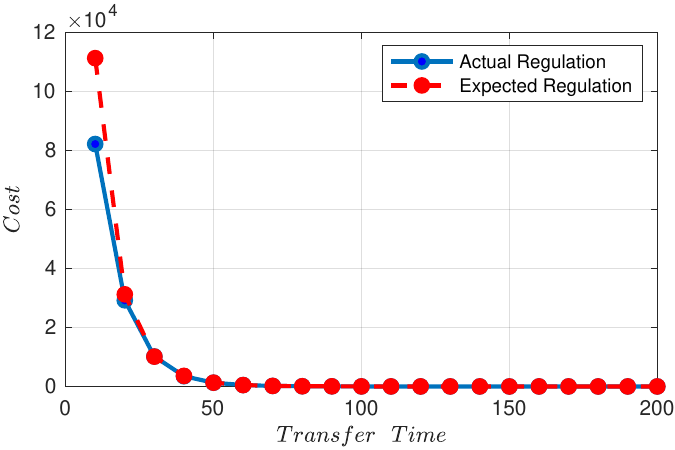}} 
    \caption{Change in total cost and regulation cost as a function of transfer time for attitude control.}
    \label{fig:att_con_cost}
\end{figure}

Figure \ref{fig:att_con_cost} illustrates the cost trajectory throughout the transfer process. It is evident that, after a certain period, the cost difference becomes negligible as the system states converge to the origin. This indicates that the regulation cost decreases significantly when the transfer time is increased, allowing the system to be linearized very close to the origin. The efficiency of this approach is underscored by the fact that as the transfer time extends, the incurred cost diminishes, ensuring a smooth and economical transition.

\begin{figure}[!htbp]
    \centering
    \subfigure[Plot of $\psi$ vs time.]{\includegraphics[width=0.32\textwidth]{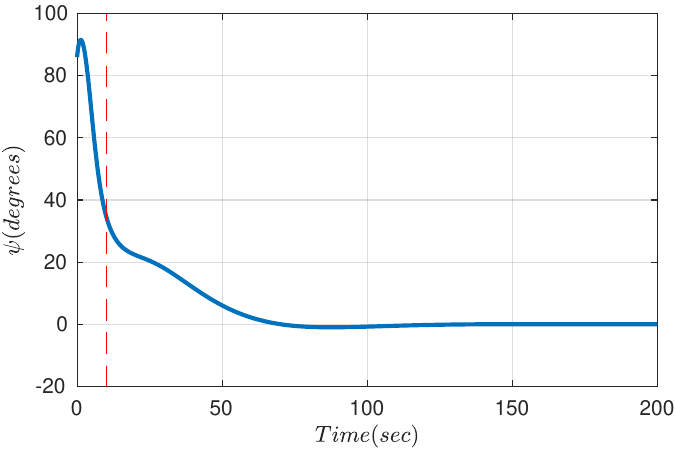}} 
    \subfigure[Plot of $\theta$ vs time.]{\includegraphics[width=0.32\textwidth]{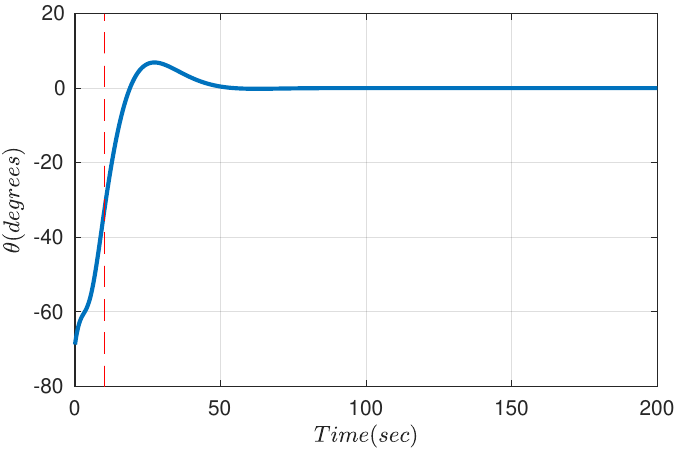}}
    \subfigure[Plot of $\phi$ vs time.]{\includegraphics[width=0.32\textwidth]{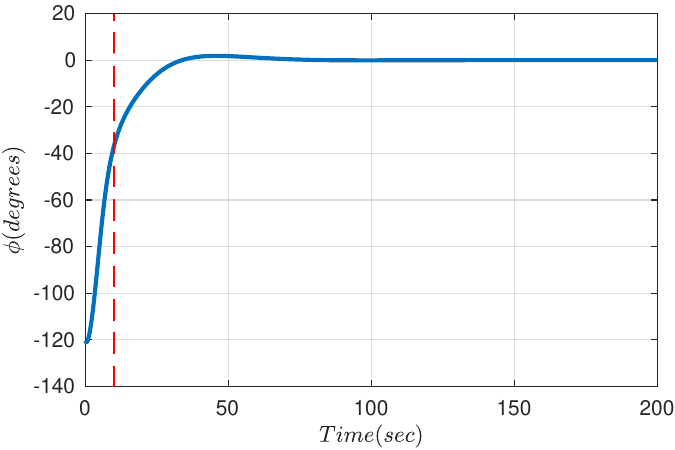}}
    \subfigure[Plot of $\omega_{1}$ vs time.]{\includegraphics[width=0.32\textwidth]{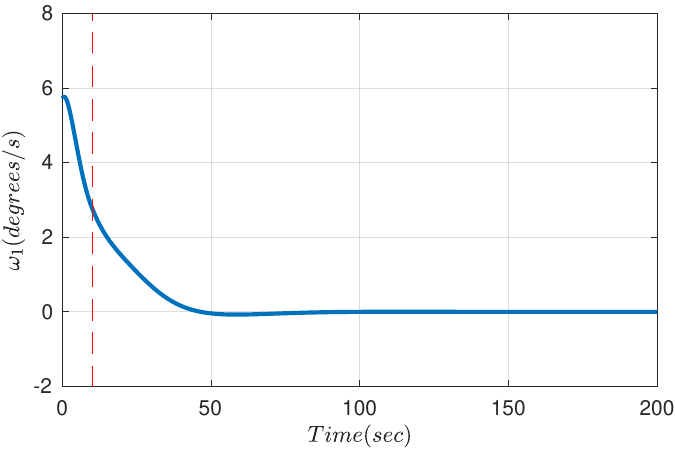}} 
    \subfigure[Plot of $\omega_{2}$ vs time.]{\includegraphics[width=0.32\textwidth]{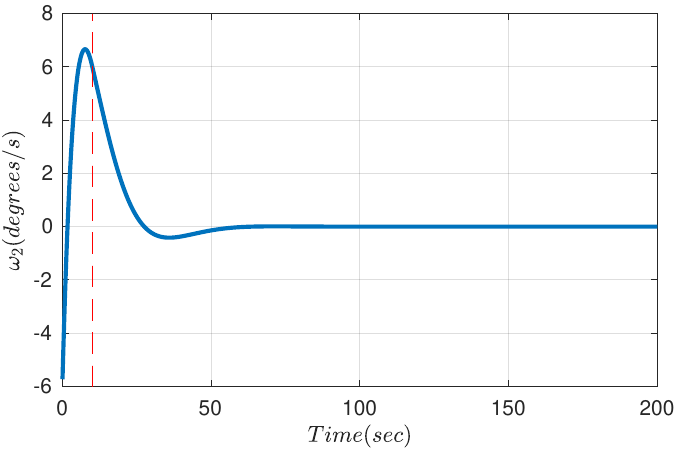}}
    \subfigure[Plot of $\omega_{3}$ vs time.]{\includegraphics[width=0.32\textwidth]{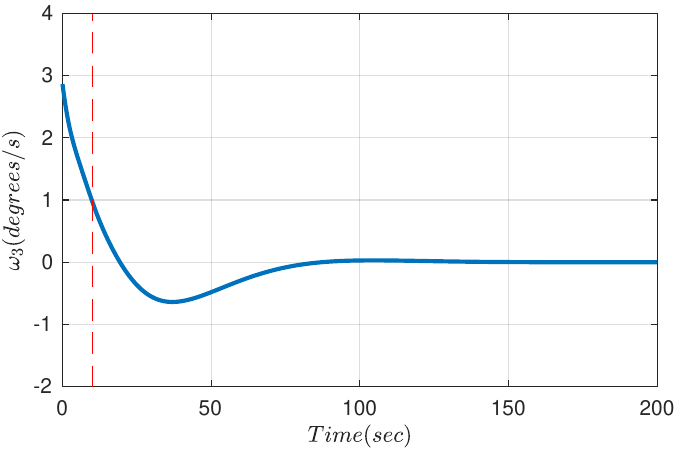}}
    \subfigure[Plot of control input $M_{1}$ vs time.]{\includegraphics[width=0.32\textwidth]{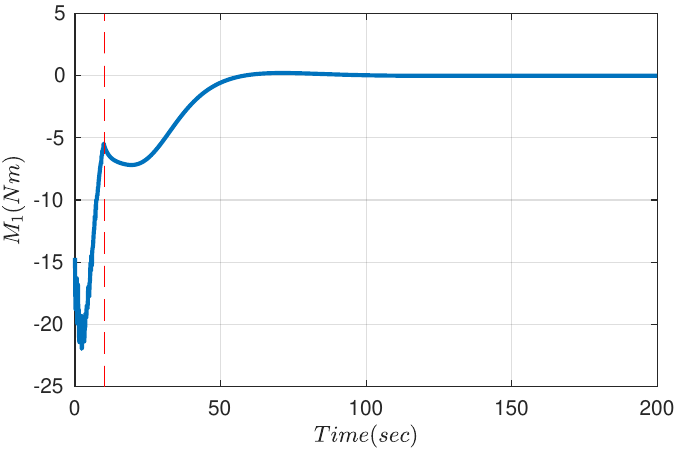}} 
    \subfigure[Plot of control input $M_{2}$ vs time.]{\includegraphics[width=0.32\textwidth]{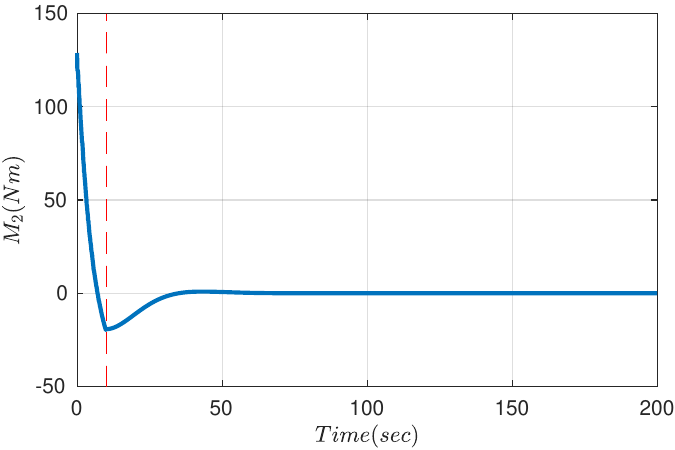}}
    \subfigure[Plot of control input $M_{3}$ vs time.]{\includegraphics[width=0.32\textwidth]{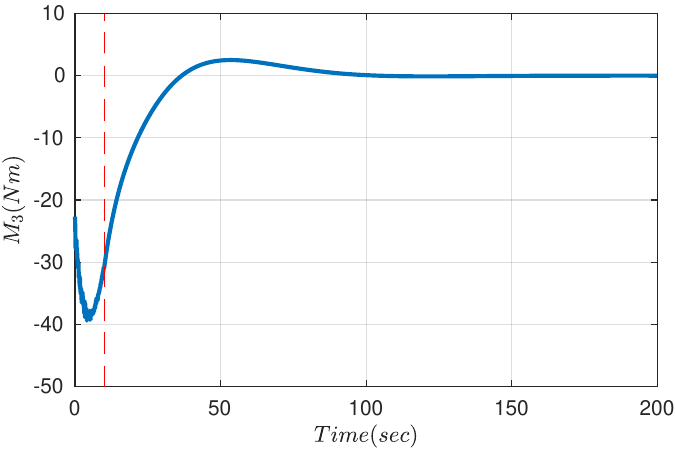}}
    \caption{Evolution of states and control for the attitude control problem with a transfer time of 10 seconds. The dotted red line shows the transition in the above plots.}
    \label{fig:att_con_states_10}
\end{figure}

\begin{figure}[!htbp]
    \centering
    \subfigure[Plot of $\psi$ vs time.]{\includegraphics[width=0.32\textwidth]{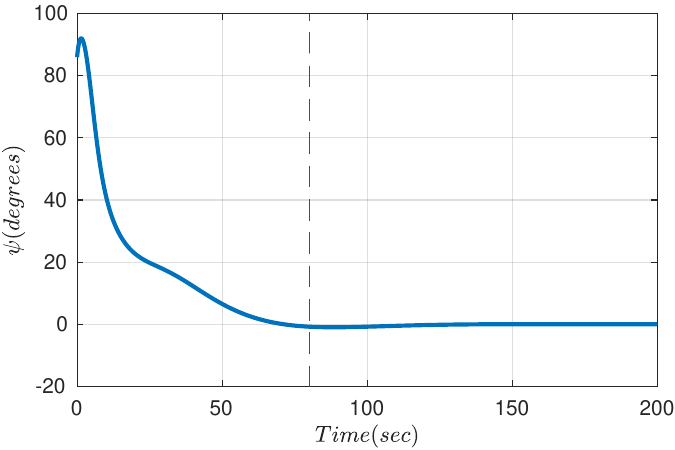}} 
    \subfigure[Plot of $\theta$ vs time.]{\includegraphics[width=0.32\textwidth]{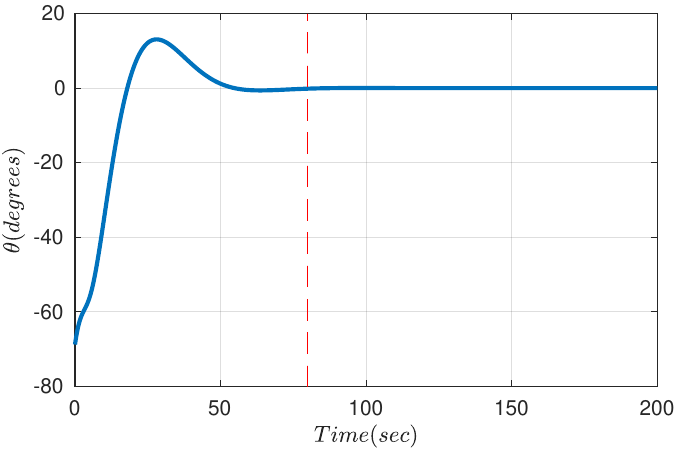}}
    \subfigure[Plot of $\phi$ vs time.]{\includegraphics[width=0.32\textwidth]{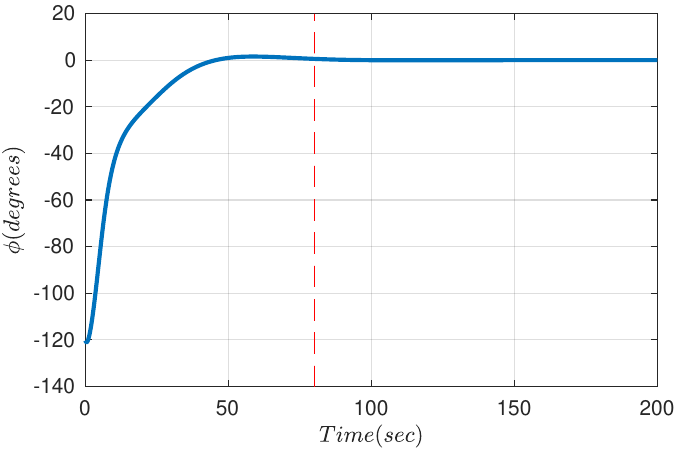}}
    \subfigure[Plot of $\omega_{1}$ vs time.]{\includegraphics[width=0.32\textwidth]{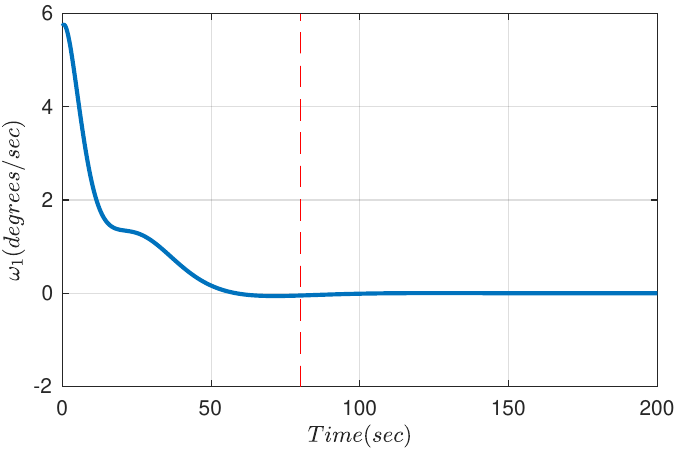}} 
    \subfigure[Plot of $\omega_{2}$ vs time.]{\includegraphics[width=0.32\textwidth]{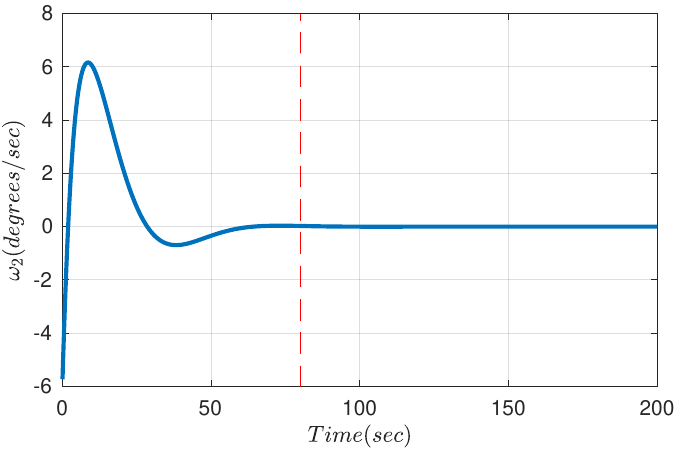}}
    \subfigure[Plot of $\omega_{3}$ vs time.]{\includegraphics[width=0.32\textwidth]{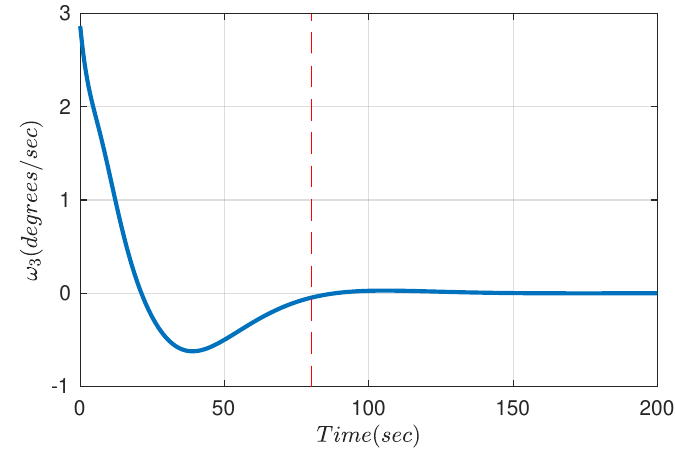}}
    \subfigure[Plot of control input $M_{1}$ vs time.]{\includegraphics[width=0.32\textwidth]{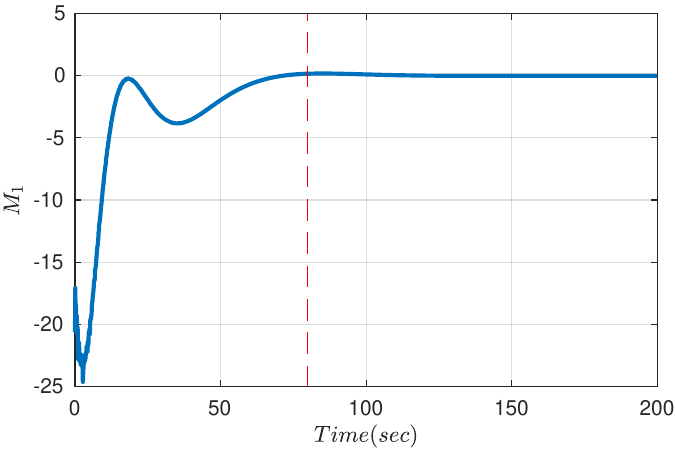}} 
    \subfigure[Plot of control input $M_{2}$ vs time.]{\includegraphics[width=0.32\textwidth]{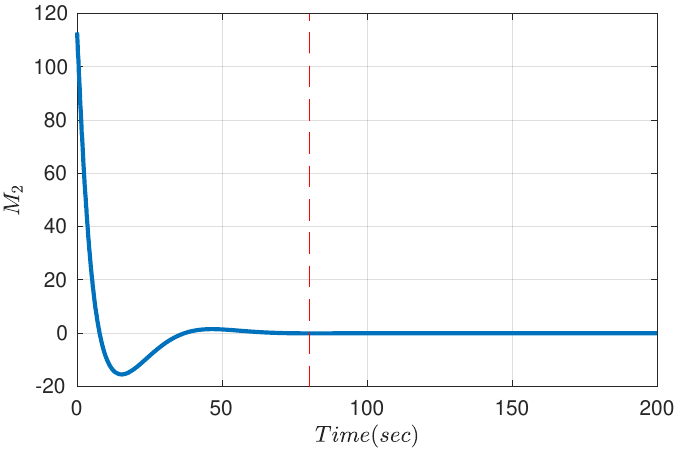}}
    \subfigure[Plot of control input $M_{3}$ vs time.]{\includegraphics[width=0.32\textwidth]{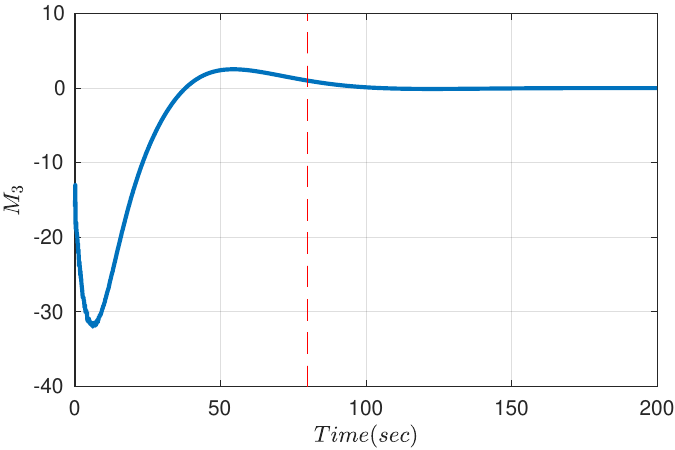}}
    \caption{Evolution of states and control for the attitude control problem with a transfer time of 80 seconds. The dotted red line shows the transition in the above plots.}
    \label{fig:att_con_states_80}
\end{figure}

To provide a deeper insight, two specific simulations are presented in Figures \ref{fig:att_con_states_10} and \ref{fig:att_con_states_80}. These figures demonstrate the impact of linearizing the system at different distances from the origin. In Figure \ref{fig:att_con_states_10}, the transfer time is set to 10 seconds and the final state errors from iLQR is [$34.8609^{o}$, $-33.1920^{o}$, $-36.7104^{o}$, $2.7864^{o}/s$, $6.0248^{o}/s$, $0.9728^{o}/s$]. Here, the states remain far from the origin at the transfer, resulting in a significantly higher cost of $6.45 \times 10^5$. This is contrasted with Figure \ref{fig:att_con_states_80}, where a transfer time of 80 seconds is considered, reducing the cost to $5.20 \times 10^5$. The final state errors from iLQR are [$-0.7695^{o}$, $-0.1469^{o}$, $0.5468^{o}$, $-0.0488^{o}/s$, $0.0243^{o}/s$, $-0.0485^{o}/s$]. This comparison clearly shows that longer transfer times lead to reduced costs and smoother transitions.

A closer examination of Figure \ref{fig:att_con_states_10} reveals a sharp change in the control input $M_{1}$. This abrupt change is indicative of the high cost and inefficiency associated with linearizing the system too far from the origin. When the transfer time is only 10 seconds, the point of linearization results in a substantial error in the state variables. Conversely, for an 80-second transfer time, this error is minimized to a very small value, leading to smoother and more cost-effective control inputs.

These observations underscore the importance of selecting an appropriate transfer time and linearization point. By increasing the transfer time, the system can be linearized closer to the origin, significantly reducing the regulation cost and ensuring smoother control transitions. This strategy not only minimizes the cost but also enhances the overall stability and efficiency of the spacecraft's maneuvering process. This detailed analysis demonstrates the critical balance between transfer time and control optimization for achieving optimal performance in attitude control.

\subsection{Rendezvous Maneuver}
For this problem, the system dynamics are defined by eq. \eqref{eq:rend}. The parameter $\alpha$ is chosen as $5 \times 10^{-4}$. Table \ref{tab:orb_param} lists the initial orbital parameters of both the chaser and target spacecraft. The value of specific impulse is chosen as $225s$ and $g_{ref}$ is $3.7114$. The time horizon for this problem is selected to be $6000$ seconds, and the system is discretized using a time step of $2$ seconds.

\begin{figure}
    \centering
    \subfigure[Total cost incurred as a function of transfer time.]{\includegraphics[width=0.48\textwidth]{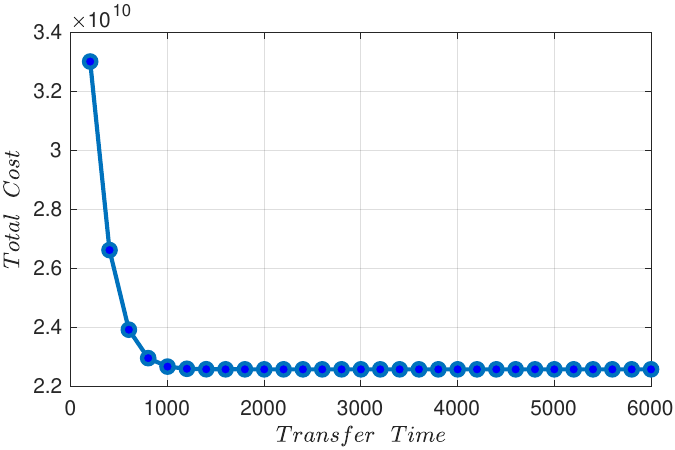}} 
    \subfigure[Cost incurred by linearized system as a function of transfer time.]{\includegraphics[width=0.48\textwidth]{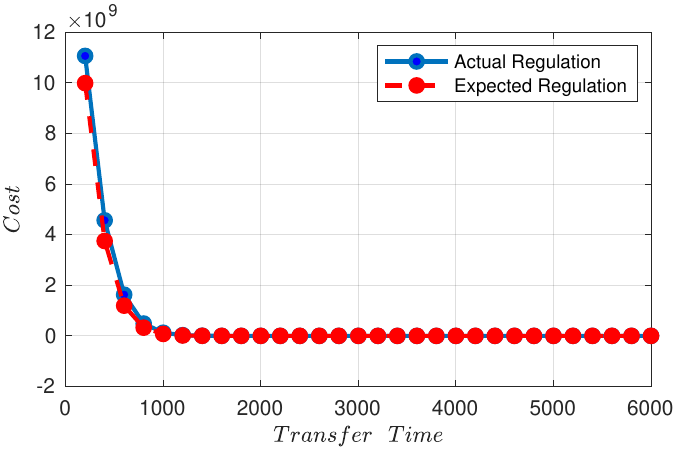}} 
    \caption{Change in total cost and regulation cost as a function of transfer time for rendezvous manuever.}
    \label{fig:rend_cost}
\end{figure}

\begin{figure}[!htbp]
    \centering
    \subfigure[Error in radial distance, $e_{r_{1}}$, vs time.]{\includegraphics[width=0.32\textwidth]{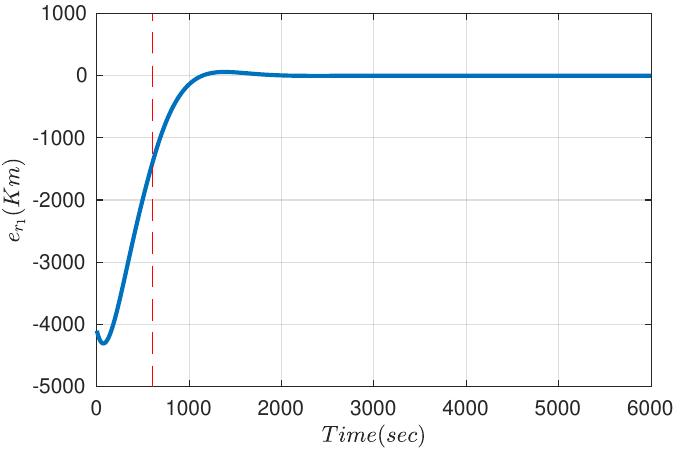}} 
    \subfigure[Error in radial distance, $e_{r_{2}}$, vs time.]{\includegraphics[width=0.32\textwidth]{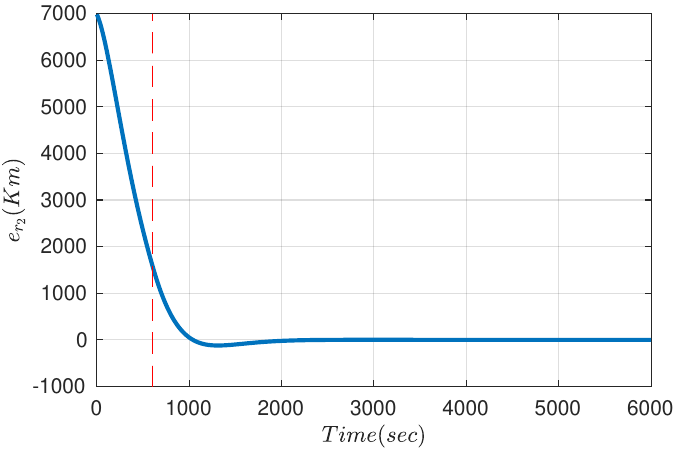}}
    \subfigure[Error in radial distance, $e_{r_{3}}$, vs time.]{\includegraphics[width=0.32\textwidth]{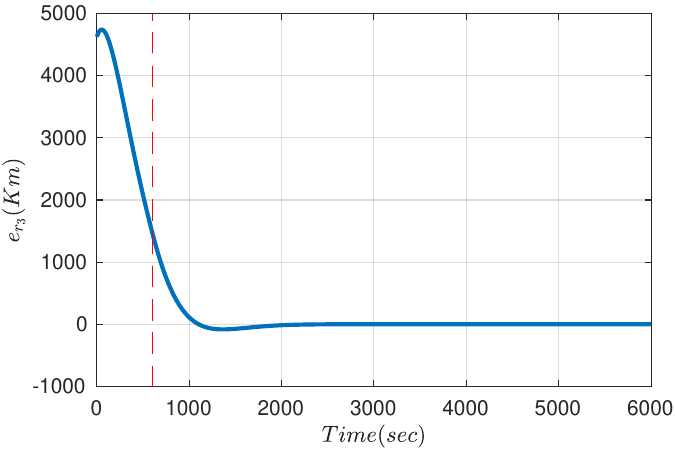}}
    \subfigure[Error in velocity, $e_{v_{1}}$, vs time.]{\includegraphics[width=0.32\textwidth]{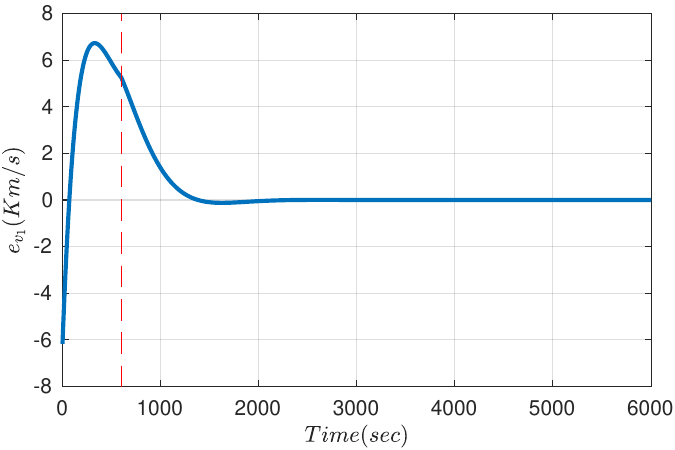}} 
    \subfigure[Error in velocity, $e_{v_{2}}$, vs time.]{\includegraphics[width=0.32\textwidth]{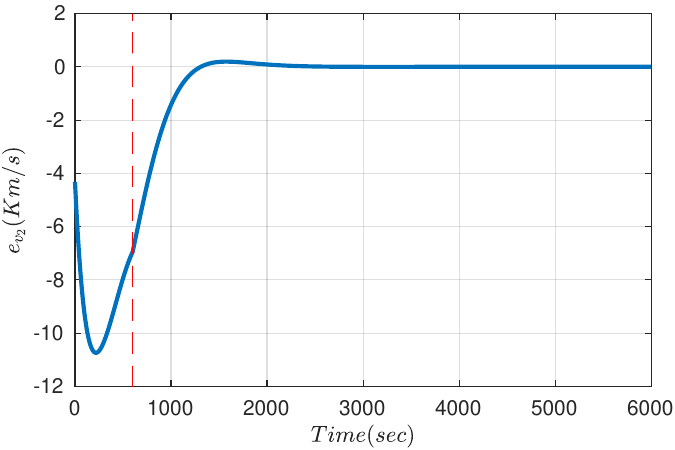}}
    \subfigure[Error in velocity, $e_{v_{3}}$, vs time.]{\includegraphics[width=0.32\textwidth]{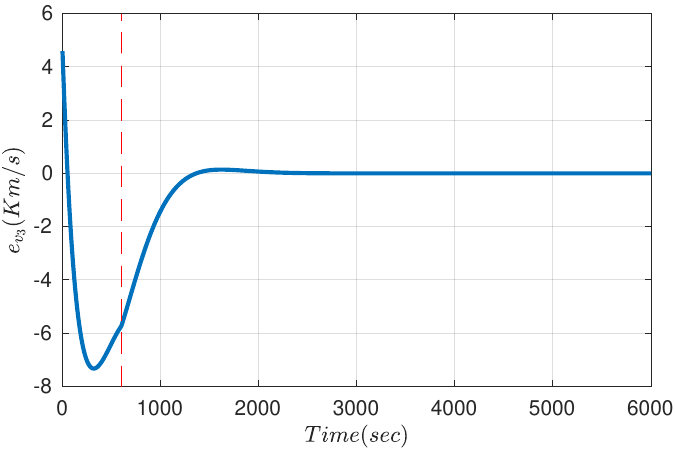}}
    \subfigure[Thrust, $u_{1}$, vs time.]{\includegraphics[width=0.32\textwidth]{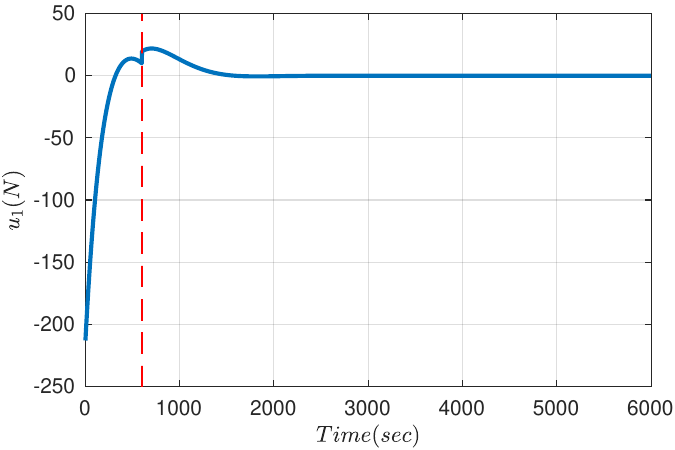}} 
    \subfigure[Thrust, $u_{2}$, vs time.]{\includegraphics[width=0.32\textwidth]{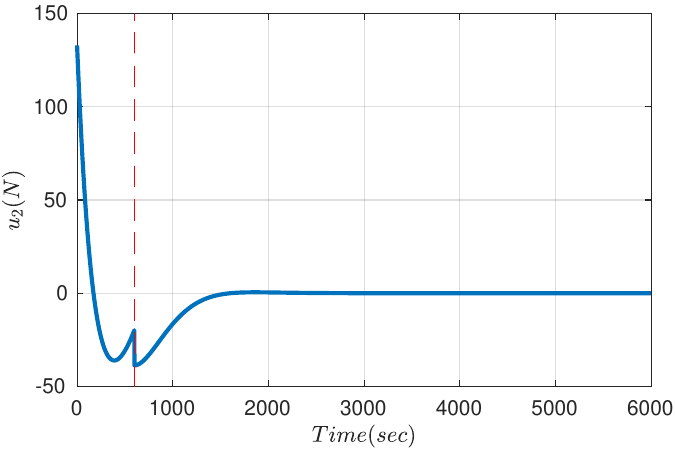}}
    \subfigure[Thrust, $u_{3}$, vs time.]{\includegraphics[width=0.32\textwidth]{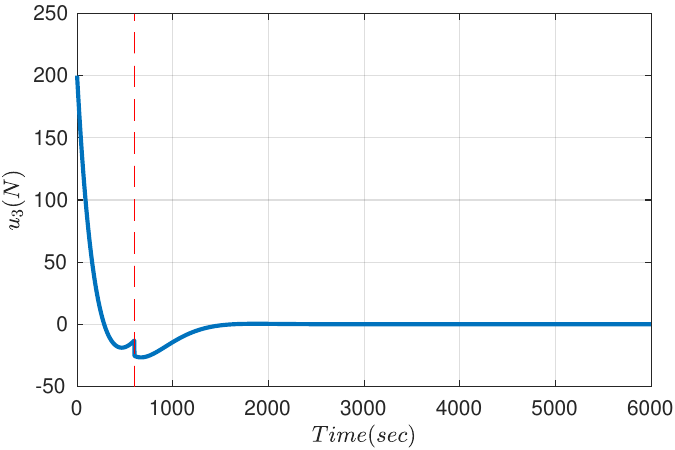}}
    \caption{Evolution of states and control for the attitude control problem with a transfer time of 600 seconds. The dotted red line shows the transition in the above plots.}
    \label{fig:rend_states_600}
\end{figure}

\begin{figure}[!htbp]
    \centering
    \subfigure[Error in radial distance, $e_{r_{1}}$, vs time.]{\includegraphics[width=0.32\textwidth]{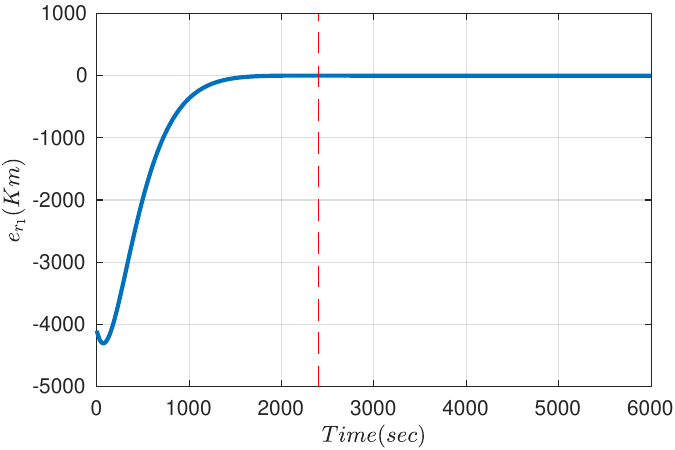}} 
    \subfigure[Error in radial distance, $e_{r_{2}}$, vs time.]{\includegraphics[width=0.32\textwidth]{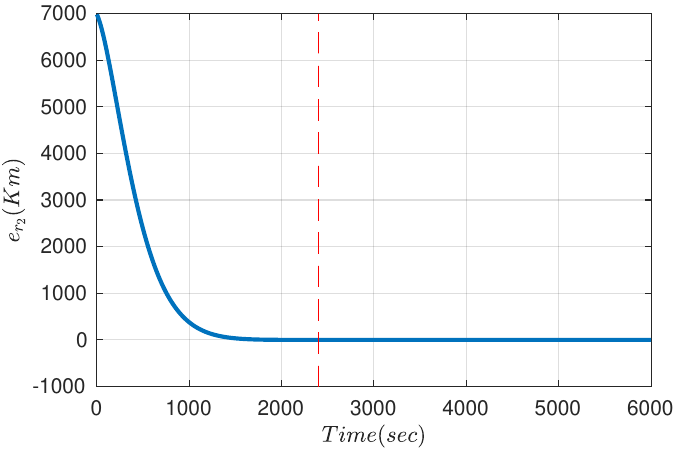}}
    \subfigure[Error in radial distance, $e_{r_{3}}$, vs time.]{\includegraphics[width=0.32\textwidth]{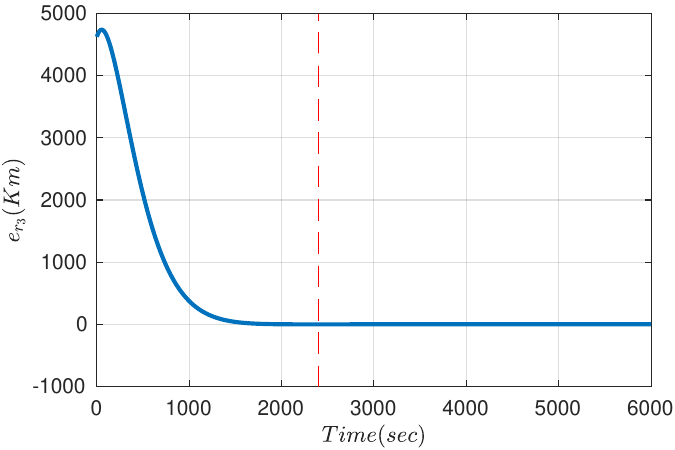}}
    \subfigure[Error in velocity, $e_{v_{1}}$, vs time.]{\includegraphics[width=0.32\textwidth]{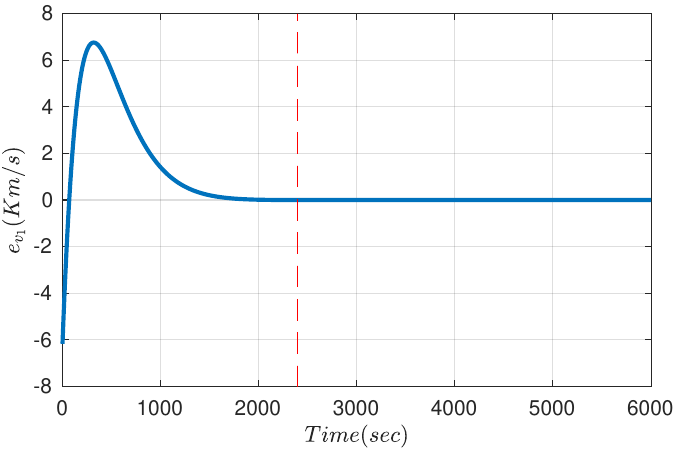}} 
    \subfigure[Error in velocity, $e_{v_{2}}$, vs time.]{\includegraphics[width=0.32\textwidth]{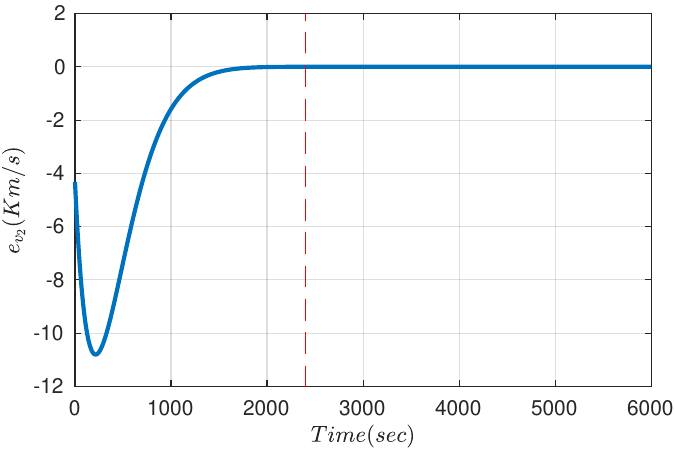}}
    \subfigure[Error in velocity, $e_{v_{3}}$, vs time.]{\includegraphics[width=0.32\textwidth]{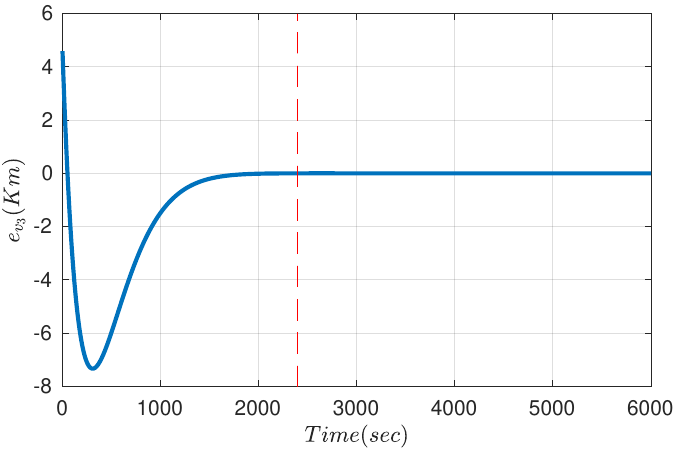}}
    \subfigure[Thrust, $u_{1}$, vs time.]{\includegraphics[width=0.32\textwidth]{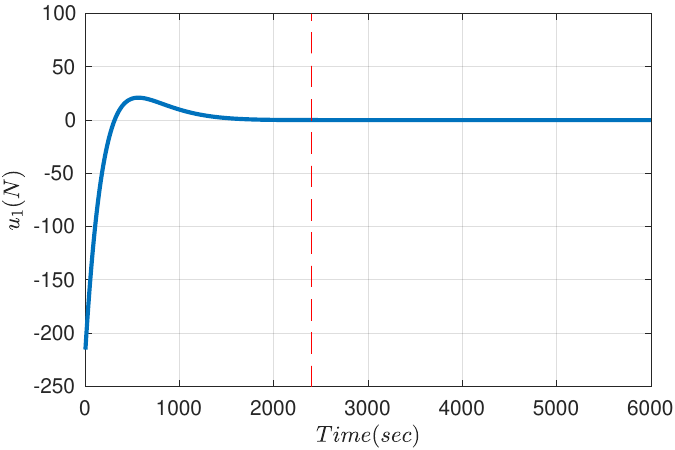}} 
    \subfigure[Thrust, $u_{2}$, vs time.]{\includegraphics[width=0.32\textwidth]{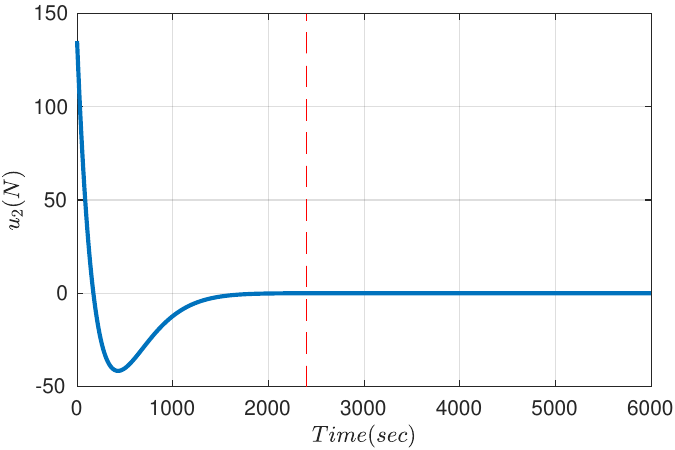}}
    \subfigure[Thrust, $u_{3}$, vs time.]{\includegraphics[width=0.32\textwidth]{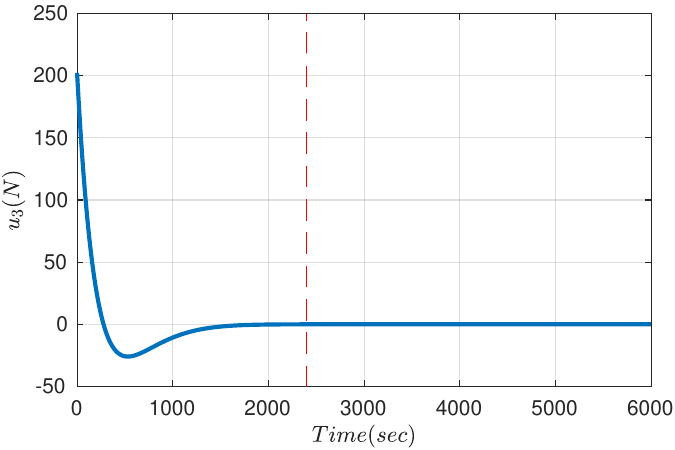}}
    \caption{Evolution of states and control for the attitude control problem with a transfer time of 2400 seconds. The dotted red line shows the transition in the above plots.}
    \label{fig:rend_states_2400}
\end{figure}

Figure \ref{fig:rend_cost} illustrates the reduction in total cost as the transfer time decreases, while the regulation cost approaches zero as the transfer time increases, indicating that the error in the final states approaches zero. To provide a comprehensive analysis, we present results for two different transfer times. Figures \ref{fig:rend_states_600} and \ref{fig:rend_states_2400} show the errors in distance and velocity, along with the control input required to drive the system, for transfer times of 600 seconds and 2400 seconds, respectively.

Similar to the behavior observed in the attitude control problem, the rendezvous maneuver demonstrates that the cost is higher for shorter transfer times, whereas it is significantly lower for longer transfer times due to reduced errors. The errors at the end of the transfer time are summarized in Table \ref{tab:Error_rend}. Despite the high error values, linear control successfully drives the system to the origin.

It is noteworthy that a sharp change in control values is observed in Figure \ref{fig:rend_states_600}, corresponding to the shorter transfer time of 600 seconds. In contrast, Figure \ref{fig:rend_states_2400} shows smoother control adjustments for the longer transfer time of 2400 seconds. This analysis underscores the importance of selecting an appropriate transfer time to minimize costs and errors in rendezvous maneuvers.

\begin{table}[!htbp]
    \centering
    \begin{tabular}{|c|c|c|c|c|c|c|}
    \hline
        Spacecraft & a & e & i & $\Omega$ & $\omega$ & $\nu$  \\
        \hline
        Chaser & 7200 Km & 0.22 & $64^{o}$ & $66^{o}$ & $28^{o}$ & $81^{o}$\\
        \hline
        Target & 7000 Km & 0.1 & $40^{o}$ & $35^{o}$ & $10^{o}$ & $120^{o}$\\
        \hline
    \end{tabular}
    \caption{Orbital parameters of chaser and target spacecraft at the initial time.}
    \label{tab:orb_param}
\end{table}

\begin{table}[!htbp]
    \centering
    \begin{tabular}{|c|c|c|c|c|c|c|}
    \hline
        Transfer Time & $e_{r_1} $(Km) & $e_{r_{2}}$ (Km)& $e_{r_{3}}$(Km) & $e_{v_{1}}$(m/s) & $e_{v_{2}}$ (m/s) & $e_{v_{3}}$ (m/s) \\
        \hline
        600 s & -1411.6 & 1601.3 & 1474.1 & 5252.7 & -6956.1 & -5732.6 \\
        \hline
        2400 s & 1.1956 & -0.8445 & -0.8217 & -2 & 1.1 & 0.5\\
        \hline
    \end{tabular}
    \caption{Error in final states at two different transfer times.}
    \label{tab:Error_rend}
\end{table}

\subsection{Soft-Landing}
The soft-landing problem, wherein we try to land a rover on the surface of Mars, poses a unique challenge when implementing our proposed approach. Computing the linear time-invariant system through linearizing about the equilibrium point becomes infeasible, since the equilibrium point includes a thrust magnitude equal to the mass times Mars's gravity. Due to this thrust, a smooth transfer to the terminal set is not possible. Moreover, this problem is constrained since $r_{3}$, representing the altitude of the lander from the surface, cannot be negative. Maintaining the lander's position requires high thrust, complicating linearization. Additionally, the lander's mass changes significantly, causing state variations over the control execution time.

Given these issues, we use the nonlinear dynamics for this problem and terminate the simulation once the altitude reaches zero. However, the altitude constraint still poses a challenge. To address this, we introduce a penalty on the altitude as a soft constraint. We add an exponential penalty function to the cost, $me^{-nr_{3}} -m$, where $r_{3}$ is the height, and $m$ and $n$ are adjustable parameters. Since the variables are of different magnitudes, we normalize them to be of the same order. We define the new variables as:
\begin{equation}
        \overline{\boldsymbol{r}} = \frac{\boldsymbol{r}}{10000},~~~
        \overline{\boldsymbol{v}} = \frac{\boldsymbol{v}}{1000},~~~
        \overline{\boldsymbol{M}} = \frac{\boldsymbol{M}}{100},~~~\text{and}~~~
        \overline{\boldsymbol{u}} = \frac{\boldsymbol{u}}{10000}.
\end{equation}
With these modifications, we reconstruct the equations and the cost function. The soft constraint introduced means the altitude might go below zero, so we terminate the simulation when the lander reaches zero altitude. For our application, we use the constraint parameters $m=100$, and $n=1$.

\begin{figure}[!htbp]
    \centering
    \subfigure[$\psi$ vs time.]{\includegraphics[width=0.32\textwidth]{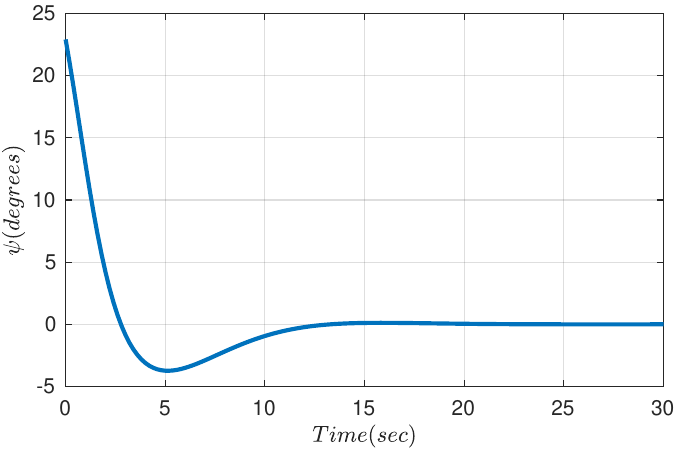}} 
    \subfigure[$\theta$ vs time.]{\includegraphics[width=0.32\textwidth]{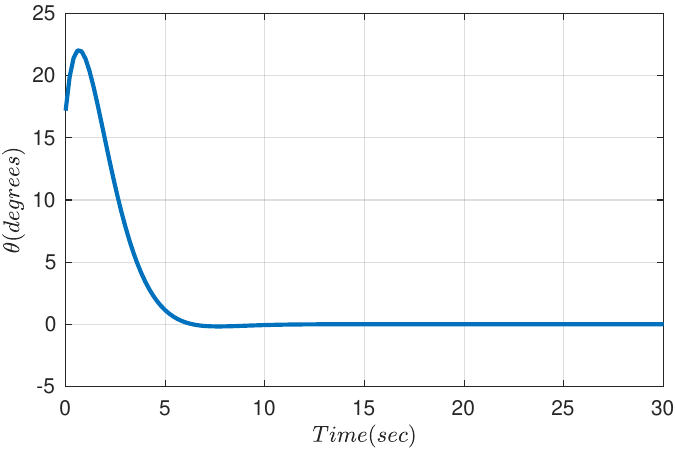}}
    \subfigure[$\phi$ vs time.]{\includegraphics[width=0.32\textwidth]{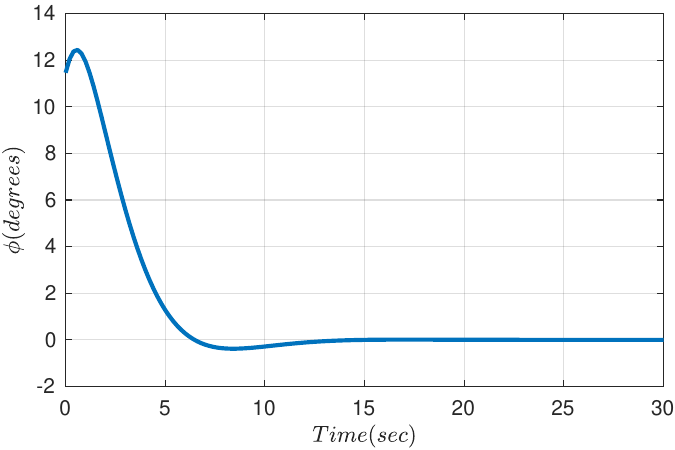}}
    \subfigure[$\omega_{1}$ vs time.]{\includegraphics[width=0.32\textwidth]{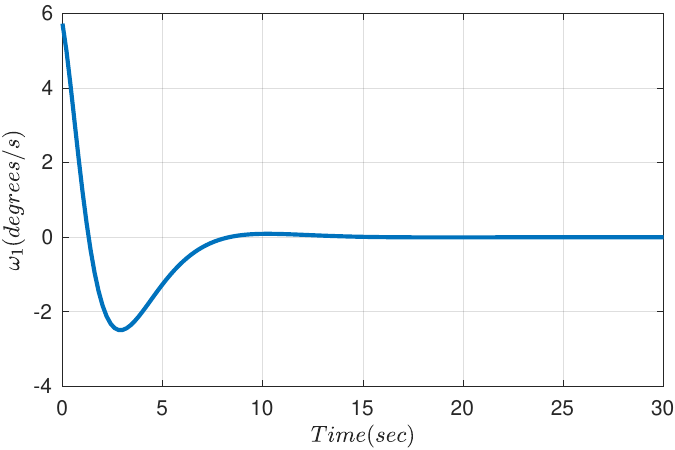}} 
    \subfigure[$\omega_{2}$ vs time.]{\includegraphics[width=0.32\textwidth]{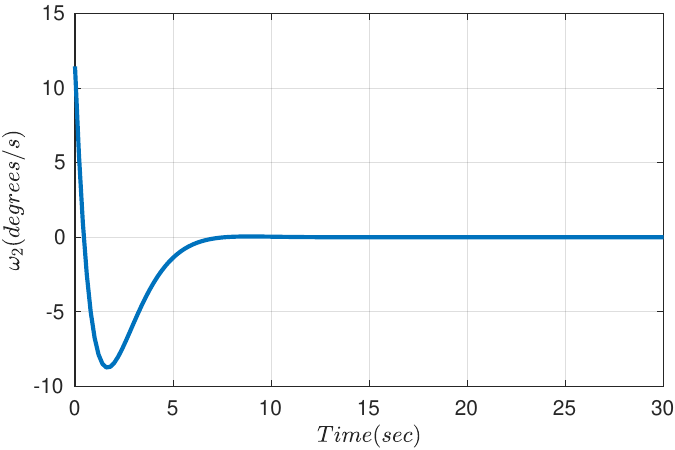}}
    \subfigure[$\omega_{3}$ vs time.]{\includegraphics[width=0.32\textwidth]{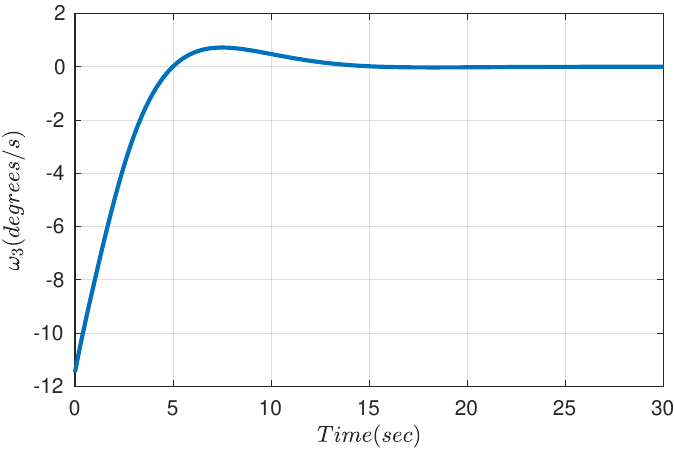}}
    \caption{Evolution of Euler angles and angular velocity for soft-landing problem.}
    \label{fig:softland_att}
\end{figure}

\begin{figure}[!htbp]
    \centering
    \subfigure[Error in position, $r_{1}$, with time.]{\includegraphics[width=0.32\textwidth]{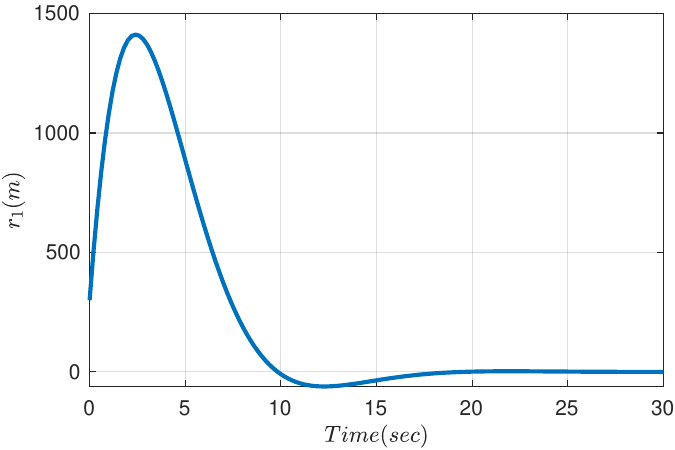}} 
    \subfigure[Error in position, $r_{2}$, with time.]{\includegraphics[width=0.32\textwidth]{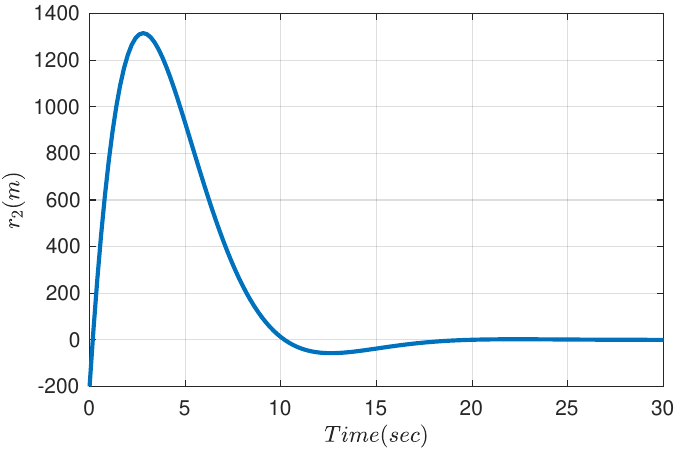}}
    \subfigure[Error in position, $r_{3}$, with time.]{\includegraphics[width=0.32\textwidth]{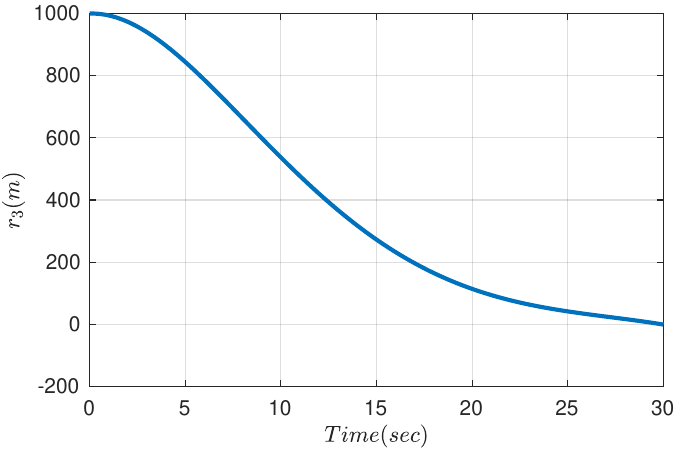}}
    \subfigure[Error in velocity, $v_{1}$, with time.]{\includegraphics[width=0.32\textwidth]{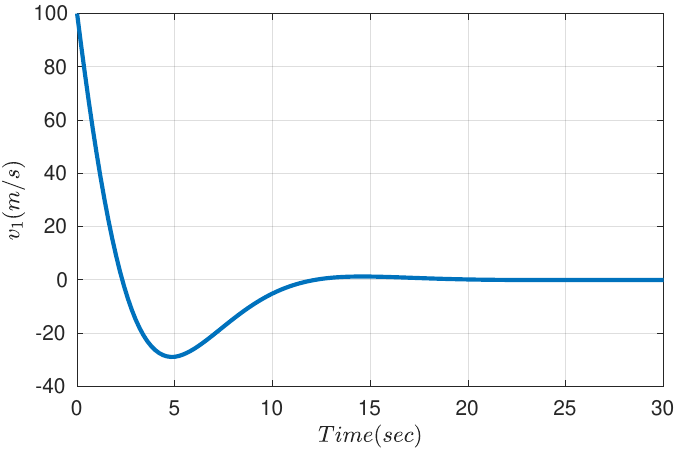}} 
    \subfigure[Error in velocity, $v_{2}$, with time.]{\includegraphics[width=0.32\textwidth]{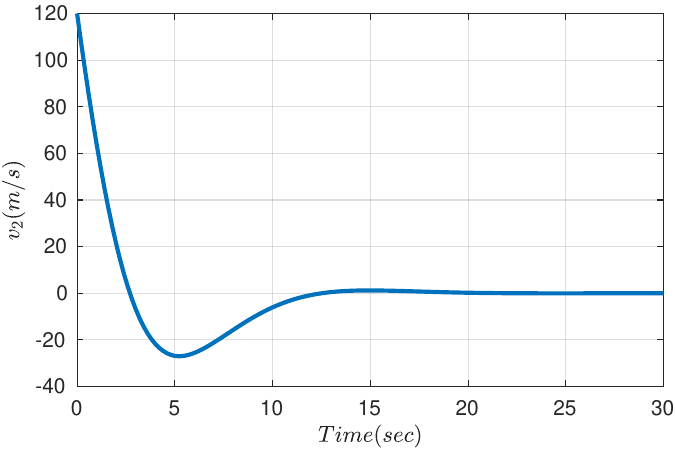}}
    \subfigure[Error in velocity, $v_{3}$, with time.]{\includegraphics[width=0.32\textwidth]{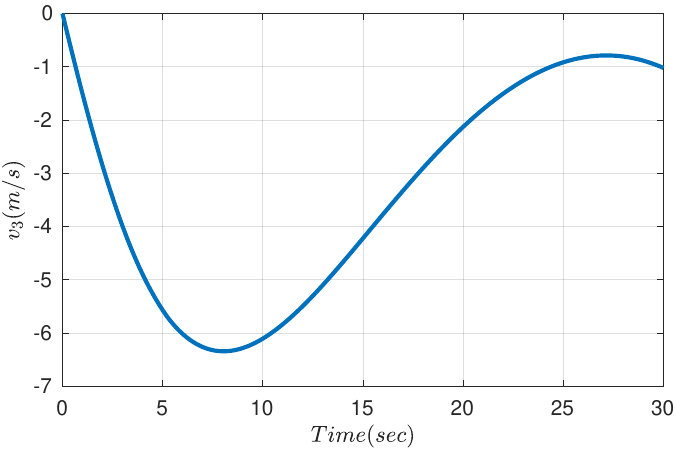}}
    \caption{Evolution of position and velocity errors for soft-landing problem.}
    \label{fig:softland_pos_vel}
\end{figure}

\begin{figure}[!htbp]
    \centering
    \subfigure[Torque input, $M_{1}$, vs time.]{\includegraphics[width=0.32\textwidth]{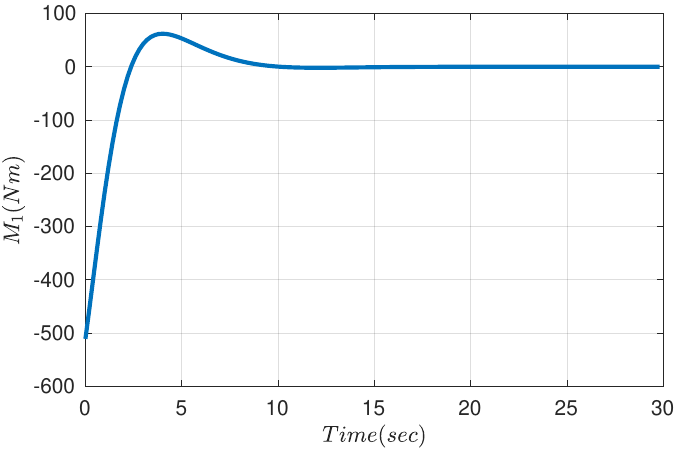}} 
    \subfigure[Torque input, $M_{2}$, vs time.]{\includegraphics[width=0.32\textwidth]{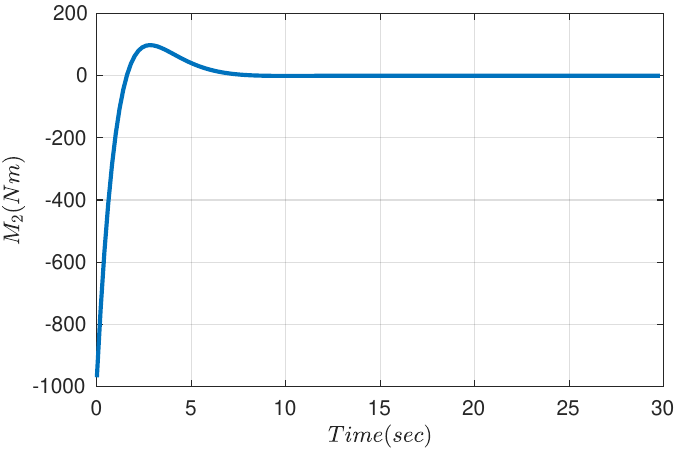}}
    \subfigure[Torque input, $M_{3}$, vs time.]{\includegraphics[width=0.32\textwidth]{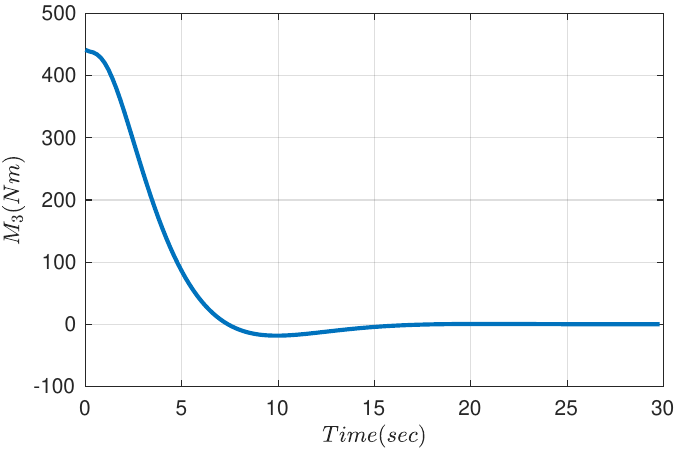}}
    \subfigure[Thrust input, $u_{1}$, vs time.]{\includegraphics[width=0.32\textwidth]{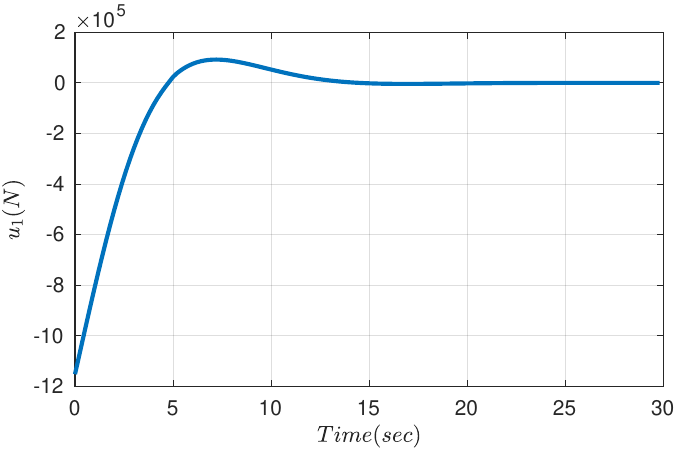}} 
    \subfigure[Thrust input, $u_{2}$, vs time.]{\includegraphics[width=0.32\textwidth]{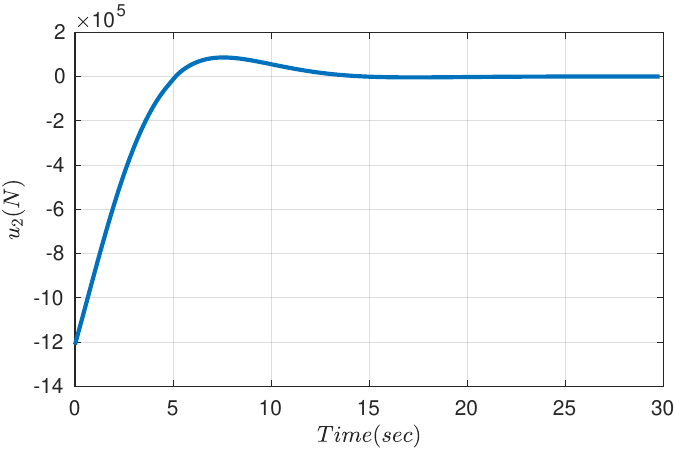}}
    \subfigure[Thrust input, $u_{3}$, vs time.]{\includegraphics[width=0.32\textwidth]{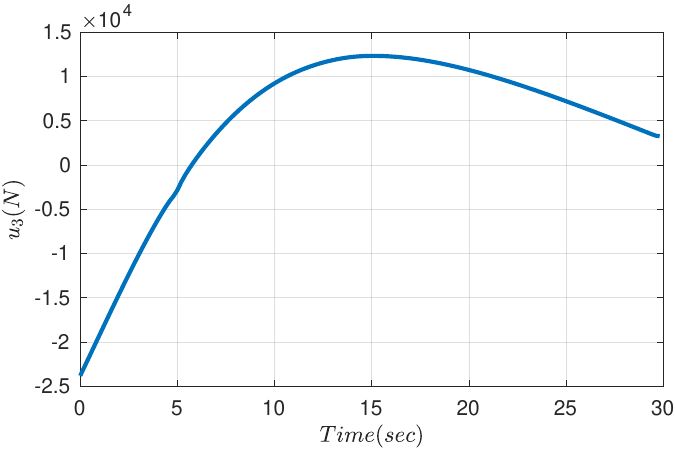}}
    \caption{Evolution of position and velocity errors for soft-landing problem.}
    \label{fig:softland_con}
\end{figure}

The initial and goal states for this problem are provided in Table \ref{tab:init_cond}, with the horizon $T=30$ seconds with a discretization step size $\delta t=0.2$ seconds. Figure \ref{fig:softland_att} shows the errors in attitude and angular velocity while soft-landing. We observe that the errors in both position and velocity converge to a safe tolerance (Fig.~\ref{fig:softland_pos_vel}). 
The final state errors in position and velocity at this point are [$-0.0595 m$, $-0.0285 m$, $1.0911 m$, $-0.0065 m/s $, $ -0.0081 m/s $, $-0.9895 m/s$]. These errors are very small for the landing to be safe. 
The lander touches the ground \textit{i.e.}, $r_{3} = 0$ at time $29.9$ seconds. So, we assume that soft-landing is completed at that point. The moment and thrust inputs are plotted in Fig.~\ref{fig:softland_con}. 

\section{Conclusion}
We pose an infinite horizon nonlinear optimal control problem and discuss a method to solve it. The algorithm capitalizes on the principle that the cost-to-go function for an infinite horizon acts as a Lyapunov function. This renders the problem globally asymptotically stable. Further, the problem is solved in two parts. We solve a finite time nonlinear OCP and then a regulation LQR problem. By transitioning a nonlinear system to a linear one near the origin, we achieve excellent performance. Our results indicate that this approach yields exceptional performance for problems such as attitude control and spacecraft rendezvous. However, for problems like soft-landing—where linearization challenges and state constraints are present—we observed that an iLQR algorithm with a penalty function proves highly effective. This method leverages state feedback to facilitate course corrections, enhancing stability. Additionally, it is crucial for the transfer to the terminal set to be smooth; otherwise, the solution may not remain optimal. This insight underscores the importance of careful transition planning in achieving consistent and reliable control performance.
\bibliographystyle{AAS_publication}   
\bibliography{naveed_references, references}   

\end{document}